\documentclass[12pt,draft,leqno]{article}
\usepackage{amssymb,eucal,latexsym}
\usepackage{amsmath,amsthm,latexsym,amscd,amsbsy,fleqn,graphicx}

\newtheorem{theorem}{Theorem}[section]
\newtheorem{lm}[theorem]{Lemma}
\newtheorem{exa}[theorem]{Example}

\newtheorem{cor}[theorem]{Corollary}
\newtheorem{pro}[theorem]{Proposition}
\newtheorem{defi}[theorem]{Definition}

\newtheorem{nota}[theorem]{Notation}

\newtheorem{rem}[theorem]{Remark}

\newtheorem{nist}[theorem]{}

\def\p{\varphi}
\def\a{\alpha}
\def\b{\beta}
\def\d{\delta}

\def\GA{\Gamma}

\def\s{\sigma}

\def\fs{\hat{f}}
\def\ps{\hat{\varphi}}

\def\lra{\longrightarrow}

\def\sbe{\subseteq}
\def\spe{\supseteq}
\def\stm{\setminus}
\def\ems{\emptyset}
\def\nes{\neq\emptyset}

\def\unl{\underline}
\def\unlx{\underline{X}}
\def\unly{\underline{Y}}
\def\unlxc{\underline{X}^c}
\def\unlb{\underline{B}}
\def\unlbc{\underline{B}^c}
\def\unla{\underline{A}}

\def\ex{\exists}
\def\fa{\forall}

\def\we{\wedge}

\def\ap{^{\prime}}
\def\inv{^{-1}}
\def\st{\ |\ }

\def\nin{\not\in}

\def\card #1{\vert #1 \vert}

\def\B{\mbox{{\boldmath $B$}}}

\def\S{{\bf S}}

\def\1{{\bf 1}}
\def\2{\mbox{{\bf 2}}}
\def\3{\mbox{{\bf 3}}}

\def\CC{{\cal C}}

\def\TT{{\cal T}}

\def\Bool{{\bf Bool}}
\def\Stone{{\bf Stone}}

\def\CCA{{\bf CCA}}
\def\CSRS{{\bf CSRS}}

\def\int{\mbox{{\rm int}}}

\def\cl{\mbox{{\rm cl}}}

\def\doc{\hspace{-1cm}{\em Proof.}~~}
\def\sq{\hspace*{\fill} \hbox{\vrule\vbox{\hrule\phantom{o}\hrule}\vrule}}
\def\sqs{\sq \vspace{2mm}}

\def\Boo{{\bf Bool}}
\def\ECS{{\bf 2Stone}}

\def\SAS{{\bf SAS}}
\def\CSAS{{\bf CSAS}}

\def\MCS{{\bf MCS}}
\def\GMCS{{\bf GMCS}}

\def\tcx{t_X^C}
\def\tcy{t_Y^C}
\def\tcx0{t_{(X,X_0)}}
\def\tcy0{t_{(Y,Y_0)}}
\def\tcxx{t_{\unlx}}
\def\tcyy{t_{\unly}}

\def\bU0{\bar{U}=(U^0,(U^i,U^{ci})_{i\in\omega})}

\def\bV0{\bar{V}=(V^0,(V^i,V^{ci})_{i\in\omega})}

\def\PCA{{\bf PCA}}
\def\PCAC{{\bf PCAC}}
\def\CA{{\bf CA}}
\def\PCS{{\bf PCS}}
\def\PCSC{{\bf PCSC}}
\def\CS{{\bf CS}}

\title{{\LARGE\bf A Generalization of}\\
\vspace{0.2cm}
{\LARGE\bf  the Stone Duality Theorem}\\
\vspace{0.5cm}
{\large\bf G. Dimov, E. Ivanova-Dimova
 and D. Vakarelov}\\
\vspace{0.2cm}
{\footnotesize\rm Department of Mathematics and Informatics,  University of Sofia,}\\
{\footnotesize\rm 5 J. Bourchier Blvd., 1164 Sofia, Bulgaria}
}

\author{}

\date{}

\begin{document}

\maketitle

\begin{abstract}
We prove a new duality theorem for the category of precontact algebras which implies the Stone Duality Theorem, its connected version obtained in \cite{DV11}, the recent duality theorems from \cite{BBSV,GG}, and some new  duality theorems for the category of contact algebras and for the category of  complete contact algebras.
\end{abstract}

\footnotetext[1]{{\footnotesize
{\em Keywords:}  (pre)contact algebra, 2-(pre)contact space,  Stone space, Stone 2-space,  (Stone) duality, C-semiregular spaces, (complete) Boolean algebra, Stone adjacency space, (closed) relations, u-points, mereocompact space.}}

\footnotetext[2]{{\footnotesize
{\em 2010 Mathematics Subject Classification:} 54E05, 18A40, 54H10, 06E15,  03G05, 54D30,  54D10.}}

\footnotetext[3]{{\footnotesize {\em E-mail addresses:}
gdimov@fmi.uni-sofia.bg, elza@fmi.uni-sofia.bg, dvak@fmi.uni-sofia.bg}}

\section{Introduction}
This paper is a continuation of the papers \cite{DV1,DV3,DV11} and, to some extent, of the papers \cite{D-APCS09, D2009,D-AMH1-10,D-AMH2-11,D2012,DI2016,DV2,VDDB}. In it we prove a new duality theorem for the category of precontact algebras which implies the Stone Duality Theorem, its connected version obtained in \cite{DV11}, the recent duality theorems from \cite{BBSV,GG}, and  some new  duality theorems for the category of contact algebras and for the category of  complete contact algebras.  More precisely, we show that there exists a duality functor $G^a$ between the category $\PCA$ of all precontact algebras and suitable morphisms between them and the category $\PCS$ of all 2-precontact spaces and suitable morphisms between them. Then, clearly, fixing some full subcategory $\cal C$ of the category $\PCA$, we obtain a duality between the categories $\CC$ and $G^a(\CC)$. Finding categories which are isomorphic or equivalent to the category $\CC$ and (or) to the category $G^a(\CC)$, we obtain as corollaries  the Stone Duality and the other dualities mentioned above. For example, when $\CC$ is the full subcategory of the category $\PCA$ having as objects all contact algebras of the form $(B,\rho_s^B)$, where $a\rho_s^B b\iff a.b\neq 0$,  we obtain that the category $\CC$ is isomorphic to the category $\Bool$ of Boolean algebras and Boolean homomorphisms and the category $G^a(\CC)$ is isomorphic to the category $\Stone$ of compact zero-dimensional Hausdorff spaces and continuous maps; in this way we obtain the Stone Duality Theorem; when $\CC$ is the full subcategory of the category $\PCA$ having as objects all contact algebras of the form $(B,\rho_l^B)$, where $a\rho_l^B b\iff a,b\neq 0$,  we obtain in a similar way the connected version of the Stone Duality Theorem, described in \cite{DV11}.

The structure of the paper is the following. In Section 2 we give the necessary preliminary results and definitions. In Section 3 we prove our main theorem, namely, the duality theorem for precontact algebras, and obtain as corollaries the theorems mentioned above.

 We now fix the notations.

 All lattices are with top (= unit) and bottom (= zero) elements,
denoted respectively by 1 and 0. We do not require the elements
$0$ and $1$ to be distinct.

 If $(X,\TT)$ is a topological space and $M$ is a subset of $X$, we
denote by $\cl_{(X,\TT)}(M)$ (or simply by $\cl(M)$ or
$\cl_X(M)$) the closure of $M$ in $(X,\TT)$ and by
$\int_{(X,\TT)}(M)$ (or briefly by $\int(M)$ or $\int_X(M)$) the
interior of $M$ in $(X,\TT)$.
The
compact spaces are not assumed to be Hausdorff (as it is adopted in \cite{E}).

If $X$ is a topological space, we denote by $CO(X)$ the set of
all clopen (=closed and open) subsets of $X$. Obviously, $(CO(X),\cup,\cap,\stm,\ems, X)$ is a Boolean algebra.

 If $X$ is a set, we denote by $2^X$ the power set of $X$.

  The set of all ultrafilters of a Boolean algebra $B$ is
denoted by $Ult(B)$.

 If $\CC$ denotes a category, we write $X\in \card\CC$ if $X$ is an
object of $\CC$, and $f\in \CC(X,Y)$ if $f$ is a morphism of $\CC$
with domain $X$ and codomain $Y$.

   The main reference books for all notions which are not defined here
 are \cite{AHS,E,kop89}.

 \section{Preliminaries}

We will first recall the notions of  {\em precontact algebra}\/ and  {\em contact algebra}.
They can be regarded as  algebraic analogues of proximity spaces (see \cite{EF,Sm,CE,AP,NW} for proximity spaces).

 \begin{defi}\label{precontact}
\rm
  An algebraic system\/ $\unlb=(B,C)$ is called a\/ {\em precontact
algebra} (\cite{DUV}) (abbreviated as PCA) if the following holds:
\begin{description}
\item[$\bullet$] \quad $B=(B,0,1,+,.,*)$ is a Boolean algebra
(where the complement is denoted by $``*$");

\item[$\bullet$]\quad  $C$ is a binary relation on $B$ (called
a {\em precontact relation}\/) satisfying the following
axioms:

\item[$(C0)$]\quad \ If $aCb$ then $a\not=0$ and $b\not=0$;

\item[$(C+)$]\quad $aC(b+c)$ iff $aCb$ or $aCc$;  $(a+b)Cc$ iff
$aCc$ or $bCc$.
\end{description}

\noindent A precontact algebra $(B,C)$ is said to be\/ {\em
complete} if the Boolean algebra $B$ is complete.
 Two precontact algebras\/ $\unlb=(B,C)$ and\/
$\underline{B_1}=(B_1,C_1)$ are said to be\/ {\em PCA-isomorphic}
(or, simply,\/ {\em isomorphic}) if there exists a {\em PCA-isomorphism} between them, i.e., a Boolean
isomorphism $\varphi: B\longrightarrow B_1$ such that, for every
$a,b\in B$, $aCb$ iff $\varphi(a)C_1 \varphi(b)$.

 The negation of the relation $C$ is denoted by $(-C)$.

 For any PCA $(B,C)$, we define a binary relation  $``\ll_C $"  on
$B$ (called {\em non-tangential inclusion})  by
\begin{equation}\label{nontangincl}
\ a \ll_C b\leftrightarrow a(-C)b^*.
\end{equation}
 Sometimes we will write simply
$``\ll$" instead of $``\ll_C$".

 We will also consider precontact algebras satisfying some
additional axioms:

\item[$(C ref)$]\quad \  If $a\not=0$ then $aCa$ (reflexivity
axiom);

\item[$(C sym)$]\quad If $aCb$ then $bCa$  (symmetry axiom);

\item[$(C tr)$]\quad\quad If $a\ll_C c$ then $(\exists b)(a\ll_C
b\ll_C c)$ (transitivity axiom);

\item[$(C con)$]\quad \  If $a\not=0,1$ then $aCa^{*}$ or $a^* Ca$
(connectedness axiom).

\smallskip

 A precontact algebra $(B,C)$ is called a\/ {\em contact algebra}
 (\cite{DV1})
(and $C$ is called a\/ {\em contact relation}) if it satisfies the
axioms $(C ref)$ and $(C sym)$.
We say that two contact algebras are {\em CA-isomorphic}\/ if they are PCA-isomorphic; also, a PCA-isomorphism between two contact algebras will be called a {\em CA-isomorphism}.

 A precontact algebra $(B,C)$ is
called\/ {\em connected}  if it satisfies the axiom $(Ccon)$.
\end{defi}

\begin{lm}\label{lemma1} Let $(B,C)$ be a precontact algebra. Define
$$aC^{\#}b\iff ((aCb)\vee  (bCa) \vee (a.b\not=0)).$$ Then $C^{\#}$ is a
contact relation on $B$ and hence $(B,C^{\#})$ is a contact
algebra.
\end{lm}

\begin{rem}\label{ctrdiez}
\rm
We will also consider precontact algebras satisfying the following
variant of the transitivity axiom (Ctr):

\smallskip

\noindent($C tr\#$)\quad\quad If $a\ll_{C^\#} c$ then $(\exists
b)(a\ll_{C^\#} b\ll_{C^\#} c)$.

\smallskip

The axiom ($C tr\#$) is known  as the $``$Interpolation axiom".

A contact algebra $(B,C)$ is called a  {\it  normal contact algebra} (\cite{deV,F}) if it satisfies the  axiom ($Ctr\#$) and the following one:

\smallskip

\noindent(C6) If $a\not= 1$ then there exists $b\not= 0$ such that
$b(-C)a$.

\smallskip

\noindent  The notion of a
normal contact algebra was introduced by Fedorchuk \cite{F} (under
the name of $``${\em Boolean $\d$-algebra}") as an equivalent expression
of the notion of a {\em compingent Boolean algebra}\/ of de Vries \cite{deV} (see its definition below). We call
such algebras $``$normal contact algebras" because they form a
subclass of the class of contact algebras and naturally arise in
normal Hausdorff spaces.

{\em The relations $C$ and $\ll$ are inter-definable.} For example,
normal contact algebras could be equivalently defined (and exactly
in this way they were introduced (under the name of {\em
compingent Boolean algebras}) by de Vries in \cite{deV}) as a pair
of a Boolean algebra $B=(B,0,1,+,.,{}^*)$ and a binary
relation $\ll$ on $B$ subject to the following axioms:

\smallskip

\noindent ($\ll$1) $a\ll b$ implies $a\leq b$;\\
($\ll$2) $0\ll 0$;\\
($\ll$3) $a\leq b\ll c\leq t$ implies $a\ll t$;\\
($\ll$4) ($a\ll b$ and $a\ll c$) implies $a\ll b.c$;\\
($\ll$5) If  $a\ll c$ then $a\ll b\ll c$  for some $b\in B$;\\
($\ll$6) If $a\neq 0$ then there exists $b\neq 0$ such that $b\ll
a$;\\
($\ll$7) $a\ll b$ implies $b^*\ll a^*$.

\smallskip

Note that if $0\neq 1$ then the axiom ($\ll$2) follows from the
axioms ($\ll$3), ($\ll$4), ($\ll$6) and ($\ll$7).

\smallskip

Obviously, contact algebras could be equivalently defined as a
pair of a Boolean algebra $B$ and a binary relation $\ll$ on $B$
subject to the  axioms ($\ll$1)-($\ll$4) and ($\ll$7); then, clearly,  the relation $\ll$
satisfies also the axioms

\smallskip

\noindent($\ll$2') $1\ll 1$;\\
($\ll$4') ($a\ll c$ and $b\ll c$) implies $(a+b)\ll c$.

 It is not difficult to see that precontact algebras could be equivalently defined as a
pair of a Boolean algebra $B$ and a binary relation $\ll$ on $B$
subject to the  axioms ($\ll$2), ($\ll$2'), ($\ll$3), ($\ll$4) and ($\ll$4').

\smallskip

It is easy to see that axiom (C6) can be stated
equivalently in the form of  ($\ll$6).
\end{rem}

Now we will give some examples of precontact and contact algebras.
We will start with the {\em extremal contact relations}.

\begin{exa}\label{extrcr}
\rm Let $B$ be a Boolean algebra. Then there exist a largest and a
smallest contact relations on $B$; the largest one, $\rho_l$ (sometimes we will write $\rho_l^B$), is
defined by $$a\rho_l b \iff (a\neq 0\mbox{ and }b\neq 0),$$ and the
smallest one, $\rho_s$ (sometimes we will write $\rho_s^B$), by $$a\rho_s b \iff a. b\neq 0.$$

Note that, for $a,b\in B$, $$a\ll_{\rho_s} b \iff a\le b;$$ hence
$a\ll_{\rho_s} a$, for any $a\in B$. Thus $(B,\rho_s)$ is a normal
contact algebra.
\end{exa}

We are now going to recall the definition of an adjacency space (\cite{Galton}, \cite{DUV}) and the fact that each adjacency space
generates canonically a precontact algebra (\cite{DUV}).

By an\/ {\em adjacency space} (see \cite{Galton} and \cite{DUV}) we mean a relational  system $(W,R)$,
where
  $W$ is a non-empty set whose elements are
called\/ {\em cells}, and $R$ is a binary relation on $W$ called
the\/ {\em adjacency relation};  the subsets of $W$ are called\/
{\em regions}.

The reflexive and symmetric closure $R^{\flat}$ of $R$ is defined
as follows:
\begin{equation}\label{rflat}
xR^{\flat}y \iff ((xRy) \vee (yRx) \vee (x=y)).
\end{equation}

 A {\em precontact relation} $C_{R}$
between the  regions of an adjacency space $(W,R)$ is defined as
follows: for every $M,N\subseteq W$,
\begin{equation}\label{CR}
MC_{R} N \mbox{ iff }(\exists x\in M)(\exists y\in N)(xRy).
\end{equation}

\begin{pro}\label{prop1}{\rm (\cite{DUV})} Let $(W,R)$ be an adjacency space
and let  \  $2^W$ be the Boolean algebra of all subsets of $W$.
Then:
\begin{description}
\item[{\rm (a)}]\quad  $(2^W, C_{R})$  is a precontact algebra;

\item[{\rm (b)}]\quad  $(2^W, C_{R})$  is a contact algebra iff $R$ is a
reflexive and symmetric relation on $W$. If $R$ is a reflexive and
symmetric relation on $W$ then $C_{R}$ coincides with
$(C_{R})^{\#}$ and $C_{R^{\flat}}$;

\item[{\rm (c)}]\quad   $C_{R}$ satisfies the axiom $(C tr)$ iff $R$ is
a transitive relation  on $W$;

\item[{\rm (d)}]\quad  $C_{R}$ satisfies the axiom $(Ccon)$ iff $R$ is a
connected relation on $W$ (which means that if $x,y\in W$ and
$x\not=y$ then there is an $R$-path from $x$ to $y$ or from $y$ to
$x$).
\end{description}
\end{pro}

Clearly, Proposition \ref{prop1}(a) implies that if $B$ is a Boolean subalgebra of the Boolean algebra $2^W$, then $(B,C_R)$ is also
a precontact algebra (here (and further on), for simplicity, we denote again by $C_R$ the restriction of the relation $C_R$ to $B$).

We  recall as well that every topological space generates canonically a contact algebra.

\begin{nist}\label{rcdef}
\rm
Let $X$ be a topological space and let $RC(X)$ be the set of all
regular closed subsets of $X$ (recall that a subset $F$ of $X$ is
said to be {\em regular closed}\/ if $F=\cl(\int(F))$). Let us equip
$RC(X)$ with the following Boolean operations and  {\em contact
relation} $C_{X}$:
\begin{description}
\item[$\bullet$] \quad $F+G=F\cup G$;

\item[$\bullet$] \quad $F^{*}=\cl(X\setminus F)$;

\item[$\bullet$] \quad  $F.G=\cl(\int(F \cap G)) (= (F^{*}\cup G^{*})^{*})$;

\item[$\bullet$] \quad $0=\emptyset$, $1=X$;

\item[$\bullet$] \quad  $FC_{X}G$ iff $F\cap G\not=\emptyset$.
\end{description}
\end{nist}

The following lemma is a well-known fact.

\begin{lm}\label{rcl} Let $X$ be a topological space. Then $$(RC(X), C_{X})=(RC(X),0,1,+,.,*,C_{X})$$ is a complete contact
algebra.
\end{lm}

The contact algebras of the type $(RC(X), C_{X})$, where $X$ is a topological space, are called {\em standard  contact algebras}.

\begin{defi}\label{ub}{\rm (\cite{DV3,DV11})}
\rm
Let\/ $\unlb=(B,C)$ be a precontact algebra and let
$U_{1},U_{2}$ be ultrafilters of $B$. We set
\begin{equation}\label{ultra}
U_{1} R_{\unlb} U_{2} \mbox{ iff }(\forall a\in U_{1})(\forall b\in
U_{2})(aCb)\ \  (\mbox{i.e., iff } U_{1}\times U_{2}\subseteq C).
\end{equation}

 The relational system $(Ult(B), R_{\unlb})$ is
called the\/ {\em canonical adjacency space of}\/  $\unlb$.

We say that $U_{1},U_{2}$ are\/ {\em connected} iff
$U_{1}(R_{\unlb})^{\flat}  U_{2}$ (see (\ref{rflat}) for the notation $R^{\flat}$).
\end{defi}

\begin{defi}\label{f}{\rm (\cite{DV3,DV11})}
\rm
 Let $X$ be a non-empty topological space and $R$ be a binary relation on $X$.
  Then the pair $(CO(X), C_{R})$ (see (\ref{CR})
for $C_R$) is a precontact algebra (by Proposition \ref{prop1}(a)), called the {\em canonical precontact
algebra of the relational system $(X,R)$.}
\end{defi}

\begin{defi}\label{tas}{\rm (\cite{DV3,DV11})}
\rm
An adjacency space $(X,R)$ is called a\/ {\em topological
adjacency space} (abbreviated as TAS) if $X$ is a topological
space and $R$ is a closed relation on $X$. When $X$ is a compact
Hausdorff zero-dimensional space (i.e., when $X$ is a\/ {\em Stone
space}), we say that the topological adjacency space $(X,R)$ is
a\/ {\em Stone adjacency space}.

Two topological adjacency spaces $(X,R)$ and $(X_1,R_1)$ are said
to be\/ {\em TAS-isomorphic} if there exists a homeomorphism
$f:X\longrightarrow X_1$ such that, for every $x,y\in X$, $xRy$
iff $f(x)R_1 f(y)$.
\end{defi}

Recall that:

\noindent(a) a topological space $X$ is called a {\em Stone space}\/ if it is a compact Hausdorff zero-dimensional space;

\noindent(b) {\em the Stone space
$S(A)$ of a Boolean algebra} $A$ is the set $X=Ult(A)$ endowed
with a topology $\TT$ having as a closed base the family
$\{s_A(a)\st a\in A\}$, where
\begin{equation}\label{sofa}
s_A(a)=\{u\in X\st a\in u\},
\end{equation}
for every $a\in A$; then $$S(A)=(X,\TT)$$ is a compact Hausdorff
zero-dimensional space (i.e., it is indeed a Stone space), $s_A(A)= CO(X)$ and {\em the Stone map}
\begin{equation}\label{stonemap}
s_A:A\lra CO(X), \ \ a\mapsto s_A(a),
\end{equation}
is a Boolean
isomorphism; also, the family $\{s_A(a)\st a\in A\}$ is an open base of $(X,\TT)$. Further, for every Stone space $X$ and   for every $x\in X$, we set
\begin{equation}\label{ux}
u_x=\{P\in CO(X)\st x\in P\}
\end{equation}
(sometimes we will also write $u_x^X$ instead of $u_x$).
Then $u_x\in Ult(CO(X))$ and the map
$$f:X\lra S(CO(X)), \ \ x\mapsto u_x,$$
is a homeomorphism.

When $\unlb=(B,C)$ is a precontact algebra, the pair $(S(B),R_{\unlb})$ is said to be {\em the canonical Stone adjacency space of} $\unlb$.

\begin{theorem}\label{th3}{\rm (\cite{DV3,DV11})}
{\rm (a)} Each precontact algebra $\unlb=(B,C)$ is isomorphic to the canonical
precontact algebra $(CO(X,\TT),C_{R_{\unlb}})$ of the Stone adjacency space
$((X,\TT),R_{\unlb})$, where $(X,\TT)=S(B)$ and for every $u,v\in X$, $uR_{\unlb}v\iff u\times v\sbe C;$
the isomorphism between them is just the Stone map $s_B:B\longrightarrow
CO(X,\TT)$.
Moreover, the relation $C$ satisfies the axiom (Cref) (resp.,
(Csym); (Ctr)) iff the relation $R_{\unlb}$ is reflexive (resp.,
symmetric; transitive).

\smallskip

\noindent{\rm (b)} There exists a bijective correspondence between the class of
all, up to PCA-isomorphism, precontact algebras and the class of
all, up to TAS-iso\-mor\-ph\-ism, Stone adjacency spaces $(X,R)$; namely, for each precontact algebra $\unlb=(B,C)$,
the PCA-isomorphism class $[\unlb]$ of $\unlb$ corresponds to the TAS-isomorphism class of the canonical Stone adjacency space $(S(B), R_{\unlb})$ of $\unlb$,
and for each Stone adjacency space $(X,R)$, the TAS-isomorphism class $[(X,R)]$ of $(X,R)$ corresponds to the PCA-isomorphism class of the canonical precontact algebra $(CO(X),C_R)$ of $(X,R)$ (see (\ref{CR}) for $C_R$).
\end{theorem}

Let us also recall the following well known statement  (see, e.g., \cite{CNG},
p.271).

 \begin{lm}\label{isombool}
Let $X$ be a dense subspace of a topological space $Y$. Then the
functions $$r:RC(Y)\lra RC(X),\  F\mapsto F\cap X,$$ and
$$e:RC(X)\lra
RC(Y),\  G\mapsto \cl_Y(G),$$ are Boolean isomorphisms between Boolean
algebras $RC(X)$ and $RC(Y)$, and $e\circ r=id_{RC(Y)}$, $r\circ
e=id_{RC(X)}$. (We will sometimes write $r_{X,Y}$ (resp., $e_{X,Y}$) instead of $r$ (resp., $e$).)
\end{lm}

 \begin{defi}\label{toppair}{\rm (\cite{DV3,DV11})}
\rm
(a) Let $X$ be a topological space and $X_{0}$ be  a dense
subspace of $X$. Then the pair
 $(X,X_{0})$ is called a\/ {\em topological pair}.

\medskip

\noindent(b) Let $(X,X_0)$ be a topological pair. Then we set
\begin{equation}\label{rcx0}
  RC(X,X_{0})=\{cl_X(A)\ |\ A\in CO(X_0)\}.
  \end{equation}
\end{defi}

\begin{lm}\label{lmrcx0}{\rm (\cite{DV3,DV11})}  Let $(X,X_0)$ be a topological pair.
Then $RC(X,X_0)\subseteq RC(X)$; the set $RC(X,X_{0})$ with the
standard Boolean operations on the regular closed subsets of $X$
is a  Boolean subalgebra of  $RC(X)$; $RC(X,X_0)$ is isomorphic to
the Boolean algebra $CO(X_{0})$; the sets $RC(X)$ and
$RC(X,X_0)$  coincide iff $X_0$ is an extremally disconnected
space. If $$C_{(X,X_0)}$$ is the restriction of the contact relation
$C_{X}$ (see Lemma \ref{rcl}) to $RC(X,X_{0})$, then $(RC(X,X_{0}),
C_{(X,X_{0})})$ is a  contact subalgebra of $(RC(X), C_{X})$.
\end{lm}

\begin{defi}\label{clandef}
\rm
Let\/ $\unlb=(B,C)$ be a precontact algebra.
A non-empty subset $\Gamma$ of $B$ is called a\/ {\em  clan} if it
satisfies the following conditions:
\begin{description}
\item[ $(Clan 1)$]\quad  $0\not\in \Gamma$;

\item[ $(Clan 2)$]\quad   If $a\in\Gamma$ and $a\leq b$ then
$b\in\Gamma$;

\item[ $(Clan 3)$]\quad  If $a+b\in \Gamma$ then  $a\in \Gamma$ or
$b\in \Gamma$;

\item[ $(Clan 4)$]\quad If $a,b\in \Gamma$ then $aC^{\#}b$.
\end{description}


The set of all clans
of a precontact
algebra\/ $\unlb$ is denoted by $Clans(\unlb)$.

Recall that a non-empty subset of a Boolean algebra $B$ is called a {\em grill}\/ if it satisfies the axioms (Clan1)-(Clan3). The set of all grills of $B$ will be denoted by $Grills(B)$.
\end{defi}

We will need the following well-known lemma   (see, e.g., \cite{Thron}):

\begin{lm}\label{grilllemma} \textsc{(Grill Lemma.)} If $F$ is a filter of a Boolean algebra $B$ and $G$ is a grill of $B$ such that $F\sbe G$ then there exists an
ultrafilter $U$ of $B$ with $F \sbe U\sbe G$.
\end{lm}

\begin{nota}\label{sigx}
\rm
Let $(X,\TT)$ be a topological space, $X_0$ be a subspace of
$X$, $x\in X$ and $B$ be a subalgebra of the Boolean algebra $(RC(X),+,.,*,\ems,X)$ defined in \ref{rcdef}. We put
\begin{equation}\label{sg}
\sigma_x^B=\{F\in B\ |\ x\in F\};\ \
\Gamma_{x,X_0}=\{F\in
CO(X_0)\ |\ x\in cl_X(F)\}.
\end{equation}

We set also

\begin{equation}\label{sg1}
\nu_x^B=\{F\in B\st x\in\int_X(F)\}.
\end{equation}

When $B=RC(X)$, we will often write simply $\s_x$ and $\nu_x$ instead of, respectively, $\s_x^B$ and $\nu_x^B$; in this case we will sometimes use the notation $\s_x^X$ and $\nu_x^X$ as well.
\end{nota}

 \begin{defi}\label{precontactspace} \textsc{(2-Precontact spaces.)}{\rm (\cite{DV3,DV11})}
\rm

\noindent(a) A triple $\unlx=(X,X_0,R)$ is  called a {\em
2-precontact space} (abbreviated as PCS) if the following
conditions are satisfied:
\begin{description}
\item[$(PCS 1)$]\quad  $(X,X_{0})$ is a topological pair and
 $X$ is a  $T_{0}$-space;

\item[$(PCS 2)$]\quad  $(X_{0},R)$ is a Stone adjacency space;

\item[$(PCS 3)$]\quad $RC(X,X_0)$ is a closed base for $X$;

\item[$(PCS 4)$]\quad For every $F,G\in CO(X_0)$, $\cl_X(F)\cap
\cl_X(G)\neq\emptyset$ implies that $F(C_R)^\# G$ (see (\ref{CR})
for $C_R$);

\item[$(PCS 5)$]\quad If $\Gamma\in Clans(CO(X_0),C_R)$ then
there exists a point $x\in X$ such that $\Gamma = \Gamma_{x,X_0}$
(see (\ref{sg}) for $\Gamma_{x,X_0}$).
\end{description}

\noindent(b) Let\/ $\unlx=(X,X_0,R)$ be a 2-precontact space. Define,
for every $F,G\in RC(X,X_0)$,
$$F\ C_{\unlx}\ G \iff ((\ex x\in
F\cap X_0) (\ex y\in G\cap X_0) (xRy)).$$
 Then the
precontact algebra
$$\unlb(\unlx)= (RC(X,X_{0}), C_{\unlx})$$
 is
said to be the\/ {\em canonical precontact algebra of\/
$\unlx$}.

\smallskip

\noindent(c) A 2-precontact space\/ $\unlx=(X,X_0,R)$ is called\/ {\em
reflexive} (resp.,\/ {\em symmetric}; {\em transitive}) if the
relation $R$ is reflexive (resp., symmetric; transitive);\/
$\unlx$ is called\/ {\em connected} if the space $X$ is
connected.

\smallskip

\noindent(d)  Let\/ $\unlx=(X,X_{0},R)$ and\/
$\widehat{\unlx}=(\widehat{X},\widehat{X}_{0},\widehat{R})$
be two 2-precontact spaces. We say that\/ $\unlx$ and\/
$\widehat{\unlx}$ are\/ {\em PCS-isomorphic} (or, simply,\/
{\em isomorphic}) if there exists a homeomorphism
$f:X\longrightarrow\widehat{X}$ such that:
\begin{description}
\item[{\rm (ISO1)}]\quad   $f(X_{0})=\widehat{X}_{0}$; and

\item[{\rm (ISO2)}]\quad $(\forall x,y\in X_{0})(xRy\leftrightarrow
f(x)\widehat{R} f(y))$.
\end{description}
\end{defi}

\begin{rem}\label{rem2pre}{\rm (\cite{DV3,DV11})}
\rm
It is very easy to see that the canonical precontact algebra of a 2-precontact space, defined in Definition \ref{precontactspace}(b),  is indeed a precontact algebra.
\end{rem}

\begin{defi}\label{canspa}{\rm (\cite{DV3,DV11})}
\rm
 Let $\unlb=(B,C)$ be a precontact algebra. We associate
 with
 $\unlb$ a 2-precontact space
$$\unlx(\unlb)=(X,X_{0},R),$$
called the\/ {\em canonical
2-precontact space of\/ $\unlb$}, as follows:
\begin{description}

\item[$\bullet$] \quad  $X=Clans(\unlb)$  and $X_{0}=Ult(B)$;

\item[$\bullet$] \quad    The topology $\TT$ on the set  $X$ is defined in
the following way:  the family $$\{g_{\unlb}(a)\
| \ a\in B\},$$ where,
for any $a\in B$,
\begin{equation}\label{fg}
g_{\unlb}(a)=\{ \Gamma\in X \ | \ a\in \Gamma\},
\end{equation}
is a closed base of $\TT$. The topology on $X_0$ is the subspace topology induced by $(X,\TT)$.

\item[$\bullet$] \quad  $R=R_{\unlb}$  (see (\ref{ultra}) for the notation $R_{\unlb}$).
\end{description}
\end{defi}

\begin{pro}\label{proposition2}{\rm (\cite{DV3,DV11})}
  Let\/ $\unlb=(B,C)$ be a precontact algebra.
 Then the  canonical
2-precontact space \/ $\unlx(\unlb)=(X,X_{0},R)$ of\/ $\unlb$
  defined above is indeed a 2-precontact
space.
\end{pro}

\begin{theorem}\label{mth} \textsc{(Representation theorem for precontact algebras.)}{\rm (\cite{DV3,DV11})}
\begin{description}
\item[{\rm (a)}]\quad  Let $\unlb=(B,C)$ be a precontact algebra
and let
$\unlx(\unlb)=(X,X_0,R)$
 be the canonical 2-precontact
space of $\unlb$. Then the function $g_{\unlb}:(B,C)\lra 2^X$, defined in
(\ref{fg}), is a PCA-isomorphism  from $(B,C)$ onto the canonical
precontact algebra $(RC(X,X_0),C_{\unlx(\unlb)})$ of $\unlx(\unlb)$. The
same function $g_{\unlb}$ is a PCA-isomorphism between  contact
algebras $(B,C^\#)$ and $(RC(X,X_0),C_{(X,X_0)})$ (see Lemma \ref{lmrcx0}(a) for $C_{(X,X_0)}$). The sets
 $RC(X)$ and  $RC(X,X_0)$ coincide iff the precontact algebra
$\unlb$ is complete. The algebra $\unlb$ satisfies the
axiom $(Cref)$ (resp., $(Csym); (Ctr))$ iff the 2-precontact
space $\unlx(\unlb)$ is reflexive (resp., symmetric;
transitive). The algebra $\unlb$ is connected iff
$\unlx(\unlb)$ is connected.

\item[{\rm (b)}]\quad There exists a bijective correspondence  between the class of
all, up to PCA-isomor\-ph\-ism, (connected) precontact algebras
and  the class of all, up to PCS-isomorphism, (connected) 2-precontact spaces; namely, for every precontact algebra $\unlb$, the PCA-isomorphism class $[\unlb]$ of $\unlb$ corresponds to the PCS-isomor\-phism class $[\unlx(\unlb)]$ of the canonical 2-precontact space $\unlx(\unlb)$ of $\unlb$, and for every 2-precontact space $\unlx$, the PCS-isomorphism class $[\unlx]$ of $\unlx$ corresponds to the PCA-isomorphism class $[\unlb(\unlx)]$ of the canonical precontact algebra $\unlb(\unlx)$ of $\unlx$.
\end{description}
\end{theorem}

\begin{cor}\label{cormth}{\rm (\cite{DV11})}
If $\unlx=(X,X_0,R)$ is a 2-precontact space then $X$ is a compact space.
\end{cor}

\begin{pro}\label{subrel}{\rm (\cite{DV3,DV11})}
Let $X_0$ be a subspace of a topological space  $X$. For every
$F,G\in CO(X_0)$, set
\begin{equation}\label{deltarel}
F\delta_{(X,X_0)}G \mbox{ iff }cl_X(F)\cap cl_X(G)\neq\emptyset.
\end{equation}
Then $(CO(X_0),\delta_{(X,X_0)})$ is a contact algebra.
\end{pro}

\begin{defi}\label{contactspace}\textsc{(2-Contact spaces.)}{\rm (\cite{DV3,DV11})}
\rm
(a) A topological pair  $(X,X_0)$ is  called a\/ {\em 2-contact
space} (abbreviated as CS) if the following conditions are
satisfied:
\begin{description}
\item[$(CS 1)$]\quad
 $X$ is a  $T_{0}$-space;

\item[$(CS 2)$]\quad  $X_{0}$ is a Stone space;

\item[$(CS 3)$]\quad $RC(X,X_0)$ is a closed base for $X$;

\item[$(CS 4)$]\quad If $\Gamma\in
Clans(CO(X_0),\delta_{(X,X_0)})$ (see (\ref{deltarel}) for the notation $\delta_{(X,X_0)}$) then there exists a point
$x\in X$ such that $\Gamma =\Gamma_{x,X_0}$ (see (\ref{sg}) for
$\Gamma_{x,X_0}$).
\end{description}

\noindent A 2-contact space $(X,X_0)$ is called\/ {\em connected}
if the space $X$ is connected.

\medskip

\noindent(b) Let  $(X,X_0)$ be a 2-contact space. Then the contact algebra
 $${\unlbc(X,X_0)}= (RC(X,X_{0}), C_{(X,X_0)})$$
(see Lemma \ref{lmrcx0}(a) for the notation $C_{(X,X_0)}$)  is said to be
the\/ {\em canonical contact algebra of the 2-contact space
$(X,X_0)$}.

\medskip

\noindent(c) Let $\unlb=(B,C)$ be a contact algebra, $X=Clans(B,C)$, $X_0=Ult(B)$
and $\TT$ be the topology on $X$ described in Definition
\ref{canspa}. Take  the subspace topology on $X_0$. Then the pair
$$\unlxc(\unlb)=(X,X_0)$$
is called the\/ {\em canonical 2-contact space of the
contact algebra} $(B,C)$. (Note that $\unlxc(\unlb)$ is indeed a 2-contact space (see \cite{DV11}).)

\medskip

\noindent(d)  Let   $(X,X_{0})$ and $(\widehat{X},\widehat{X}_{0})$ be two
2-contact spaces. We say that  $(X,X_0)$ and
$(\widehat{X},\widehat{X}_0)$ are\/ {\em CS-isomorphic} (or,
simply,\/ {\em isomorphic}) if there exists a homeomorphism
$f:X\longrightarrow\widehat{X}$ such that
  $f(X_{0})=\widehat{X}_0$.
\end{defi}

\begin{cor}\label{cormthccc}{\rm (\cite{DV11})}
If $\unlx=(X,X_0,R)$ is a 2-contact space then $X$ is a compact space.
\end{cor}

\begin{defi}\label{2grillsp}\textsc{(Stone 2-spaces.)}{\rm (\cite{DV11})}
\rm

\smallskip

\noindent(a) A topological pair  $(X,X_0)$ is  called a\/ {\em Stone
2-space} (abbreviated as S2S) if it satisfies conditions (CS1)-(CS3) of Definition \ref{contactspace} and the following condition:

\begin{description}
\item[$(S2S4)$]\quad If $\Gamma\in
Grills(CO(X_0))$  then there exists a point
$x\in X$ such that $\Gamma =\Gamma_{x,X_0}$ (see (\ref{sg}) for
$\Gamma_{x,X_0}$).
\end{description}

\smallskip

\noindent(b) Let   $(X,X_{0})$ and $(\widehat{X},\widehat{X}_{0})$ be two
Stone 2-spaces. We say that  $(X,X_0)$ and
$(\widehat{X},\widehat{X}_0)$ are\/ {\em S2S-isomorphic} (or,
simply,\/ {\em isomorphic}) if there exists a homeomorphism
$f:X\longrightarrow\widehat{X}$ such that
  $f(X_{0})=\widehat{X}_0$.
\end{defi}

\begin{defi}\label{defeconsp}{\rm (\cite{DV11})}
\rm
Let $(X,X_0)$ and $(X\ap,X_0\ap)$ be two Stone 2-spaces and $f:X\lra X\ap$ be a continuous map. Then $f$ is called a {\em 2-map}\/ if $f(X_0)\sbe X_0\ap$.

The category of all Stone 2-spaces and all 2-maps between them will be denoted by $\ECS$.

The category of all Boolean algebras and all Boolean homomorphisms between them will be denoted by $\Bool$.
\end{defi}

Recall that a topological space $X$ is said to be {\em semiregular} if $RC(X)$ is a closed base for $X$.

\begin{defi}\label{csemi}{\rm (\cite{DV1})}
\rm
A semiregular $T_0$-space $(X,\TT)$ is said to be \/ {\em C-semiregular}\/ if
for every clan $\Gamma$ in $(RC(X),C_X)$ there exists a point
$x\in X$ such that $\Gamma=\sigma_x$ (see (\ref{sg}) for
$\sigma_x$).
\end{defi}

\begin{pro}\label{csemicomp}{\rm (\cite[Fact 4.1]{DV1})}
Every C-semiregular space $X$ is a compact space.
\end{pro}

 \begin{defi}\label{upoint}{\rm (\cite{DV11})}
 \rm
 Let $(X,\TT)$ be a topological space and $x\in X$. The point $x$ is said to be an {\em u-point}\/ if for every $U,V\in\TT$, $x\in\cl(U)\cap \cl(V)$ implies that $x\in\cl(U\cap V)$.
 \end{defi}

\begin{theorem}\label{novocor2}{\rm (\cite{DV11})}
 If $X$ is  C-semiregular  and
 $$X_0=\{x\in X\st x \mbox{ is an u-point of }X\},$$
 then the pair $(X,X_0)$ is a 2-contact space and $X_0$ is a dense  extremally
 disconnected compact Hausdorff subspace of $X$; moreover, $X_0$ is the unique dense  extremally
 disconnected compact Hausdorff subspace of $X$.
 \end{theorem}

\begin{lm}\label{csem2cont}{\rm (\cite{DV11})}
If $(X,X_0)$ is a 2-contact space and $X_0$ is extremally disconnected, then $X$ is C-semiregular.
\end{lm}

\section{The Main Theorem and its corollaries}

The next lemma is obvious.

\begin{lm}\label{obvle}
Let $A$ be a subalgebra of a Boolean algebra $B$ and $(B,C)$ be a precontact algebra. Then $(A,C\cap A^2)$ is a precontact algebra and if $\GA\in Clans(B,C)$ then $\GA\cap A\in Clans(A,C\cap A^2)$.
\end{lm}

We will need also a lemma from \cite{DV11,DV1}.

\begin{lm}\label{pro41dv1}{\rm (\cite{DV11, DV1})}
Let $X$ be a topological space, $B$ be a subalgebra of the Boolean algebra $(RC(X),+,.,*,\ems,X)$ (defined in \ref{rcdef}), $\GA$ be a grill of $B$, $x\in X$ and $\GA\sbe\s_x^B$. Then $\nu_x^B\sbe \GA$ (see (\ref{sg1}) for $\nu_x^B$ and (\ref{sg}) for $\s_x^B$).
\end{lm}

\begin{defi}\label{precontactdualitydef}
\rm
 Let $\PCA$ be the category of all precontact algebras and all Boolean homomorphisms $\p:(B,C)\lra(B\ap,C\ap)$ between them such that, for all $a,b\in B$, $\p(a)C\ap\p(b)$ implies that $aCb$.

 Let $\PCS$ be the category of all 2-precontact spaces and all continuous maps $f:(X,X_0,R)\lra(X\ap,X_0\ap,R\ap)$ between them such that $f(X_0)\sbe X_0\ap$ and, for every $x,y\in X_0$, $xRy$ implies that $f(x)R\ap f(y)$.
\end{defi}

The proof of the next theorem is similar in some of its parts to the proof of \cite[Theorem 8.17]{DV11} but, for convenience of the reader, we give here all necessary details.

\begin{theorem}\label{maindualitytheorem} {\sc (The Main Theorem: A Duality Theorem for Precontact Algebras)}
The categories\/ $\PCA$ and\/ $\PCS$ are dually equivalent. In fact, the contravariant functor $G^a:\PCA\lra\PCS$ defined by $G^a(\unlb)=\unlx(\unlb)$ on the objects $\unlb$ of the category\/ $\PCA$, and by $G^a(\p):G^a(\unl{B\ap})\lra G^a(\unlb)$, $\GA\mapsto\p\inv(\GA)$, for every\/ $\PCA$-morphism $\p:\unlb\lra\unl{B\ap}$, is a duality functor.
\end{theorem}

\doc   We will first define two
contravariant functors $$G^a:\PCA\lra\PCS\ \ \mbox{ and }\ \
G^t:\PCS\lra\PCA.$$

Let $\unl{X}=(X,X_0,R)\in\card{\PCS}$. Define $$G^t(X,X_0,R)=\unlb(\unl{X}),$$
i.e. $G^t(X,X_0,R)$ is the canonical precontact algebra $\unlb(\unlx)=(RC(X,X_0),C_{\unlx})$ of the  2-precontact space $\unlx$
(see Definition \ref{precontactspace}(b) and recall that for every $F_1,F_2\in RC(X,X_0)$,
$$F_1 C_{\unlx} F_2\iff ((\ex x_1\in F_1\cap X_0)(\ex x_2\in F_2\cap X_0)(x_1 R x_2))).$$
 Hence $G^t(X,X_0,R)\in\card{\PCA}$.

Let $f\in\PCS((X,X_0,R),(Y,Y_0,R\ap))$. Define $$G^t(f):G^t(Y,Y_0,R\ap)\lra G^t(X,X_0,R)$$ by
the formula
\begin{equation}\label{deftetafc}
G^t(f)(\cl_Y(F\ap))=\cl_X(X_0\cap f\inv(F\ap)), \ \  \fa F\ap\in CO(Y_0).
\end{equation}
Set $\p_f=G^t(f)$, $\unlx=(X,X_0,R)$, $\unly=(Y,Y_0,R\ap)$, $C=C_{\unlx}$ and $C\ap=C_{\unly}$. We will show that $\p_f$ is a Boolean homomorphism
between the Boolean algebras $RC(Y,Y_0)$ and $RC(X,X_0)$. Clearly, for every $F\ap\in CO(Y_0)$, we have that $X_0\cap f\inv(F\ap)=(f|_{X_0})\inv(F\ap)\in CO(X_0)$. Hence $\p_f(RC(Y,Y_0))\sbe RC(X,X_0)$. Obviously,
$\p_f(\ems)=\ems$ and $\p_f(Y)=X$.
Let $F\ap,G\ap\in CO(Y_0)$. Then $\p_f(\cl_Y(F\ap)+\cl_Y(G\ap))=\p_f(\cl_Y(F\ap\cup G\ap))=\cl_X(X_0\cap f\inv(F\ap\cup G\ap))=\cl_X((X_0\cap f\inv(F\ap))\cup(X_0\cap f\inv(G\ap)))=\p_f(\cl_Y(F\ap))+
\p_f(\cl_Y(G\ap))$. Also, using Lemma \ref{isombool}, we get that  $\p_f((\cl_Y(F\ap))^*)=\p_f(\cl_Y(Y_0\stm F\ap)=\cl_X(X_0\cap f\inv(Y_0\stm F\ap))=\cl_X(X_0\cap(f\inv(Y_0)\stm f\inv(F\ap)))=\cl_X(X_0\stm(X_0\cap f\inv(F\ap)))=(\cl_X(X_0\cap f\inv(F\ap)))^*=(\p_f(\cl_Y(F\ap))^*$. So,
\begin{equation}\label{deftetafcnov}
G^t(f) \ \mbox{ is a Boolean homomorphism}.
\end{equation}

Let $F\ap,G\ap\in CO(Y_0)$ and $\p_f(\cl_Y(F\ap)) C \p_f(\cl_Y(G\ap))$. Then
$$\cl_X(X_0\cap f\inv(F\ap)) C \cl_X(X_0\cap f\inv(G\ap))$$
 and thus there exist $x\in X_0\cap\cl_X(X_0\cap f\inv(F\ap))(= X_0\cap f\inv(F\ap))$ and $y\in X_0\cap \cl_X(X_0\cap f\inv(G\ap))(= X_0\cap f\inv(G\ap))$ such that $xRy$. Hence $f(x)R\ap f(y)$. By Lemma \ref{isombool}, we get that $F\ap=Y_0\cap\cl_Y(F\ap)$ and $G\ap=Y_0\cap\cl_Y(G\ap)$. Since, obviously, $f(x)\in F\ap$ and $f(y)\in G\ap$, we get that
 $\cl_Y(F\ap)C\ap\cl_Y(G\ap)$. Thus $\p$ is a $\PCA$-morphism,
i.e., $G^t(f)$ is well defined.

Now we will show  that $G^t$ is a contravariant functor. Clearly, $G^t(id_{(X,X_0,R)})=id_{G^t(X,X_0,R)}$.

Let $f\in\PCS((X,X_0,R),(Y,Y_0,R\ap))$ and $g\in\PCS((Y,Y_0,R\ap),(Z,Z_0,R''))$.
Then, for every $F\in CO(Z_0)$, $G^t(g\circ f)(\cl_Z(F))=\cl_X(X_0\cap(g\circ f)\inv(F))=\cl_X(X_0\cap f\inv(g\inv(F)))$ and $(G^t(f)\circ G^t(g))(\cl_Z(F))=G^t(f)(\cl_Y(Y_0\cap g\inv(F)))=\cl_X(X_0\cap f\inv(Y_0\cap g\inv(F)))=\cl_X(X_0\cap f\inv(Y_0)\cap f\inv(g\inv(F)))=\cl_X(X_0\cap f\inv(g\inv(F)))=
G^t(g\circ f)(\cl_Z(F))$. So, $G^t$ is a contravariant functor.

For every precontact algebra $\unlb$, set $$G^a(\unlb)=\unlx(\unlb),$$
where $\unlx(\unlb)$ is the canonical  2-precontact space of the precontact  algebra $\unlb$ (see Definition \ref{canspa}).
Then Proposition \ref{proposition2} implies that
$G^a(\unlb)\in\card{\PCS}$.

Let $\unlb=(B,C)$,  $\unlb\ap=(B\ap,C\ap)$ and $\p\in\PCA(\unlb,\unlb\ap)$. Let $G^a(\unlb)=(X,X_0,R)$ and $G^a(\unlb\ap)=(Y,Y_0,R\ap)$.
Then we define the map
$$G^a(\p):G^a(\unlb\ap)\lra G^a(\unlb)$$ by the formula
\begin{equation}\label{deftetaphic}
G^a(\p)(\GA)=\p\inv(\GA), \ \  \fa \GA\in Y.
\end{equation}
Set $f_\p=G^a(\p)$. Since every grill of a Boolean algebra $B\ap$ is a union of ultrafilters of $B\ap$, every union of ultrafilters of $B\ap$ is a grill of $B\ap$ (see, e.g., \cite[Corollary 3.1]{DV1}), and the inverse image of an ultrafilter by a Boolean homomorphism between two Boolean algebras is again an ultrafilter, we get that
$\fa \GA\in Y$, $f_\p(\GA)$ is a grill of $B$. We have to show that $f_\p(\GA)$ is a clan. Let $a,b\in f_\p(\GA)$, i.e. $a,b\in\p\inv(\GA)$. Hence $\p(a)\in\GA\ni\p(b)$.
Thus $\p(a)(C\ap)^{\#}\p(b)$, i.e. $(\p(a)C\ap\p(b))\vee(\p(b)C\ap\p(a))\vee(\p(a).\p(b)\neq 0)$. Then $(aCb)\vee(bCa)\vee(a.b\neq 0)$. Therefore, $aC^{\#}b$. This shows that $f_\p(\GA)$ is a clan,
i.e. $f_\p:Y\lra X$.

We will show that $f_\p$ is a continuous function. Let $a\in B$. Then $g_{\unlb}(a)=\{\GA\in X\st a\in \GA\}$ is a basic closed subset of $X$ (see Definition \ref{canspa}). We will show that
\begin{equation}\label{contfphic}
f_\p\inv(g_{\unlb}(a))=g_{\unlb\ap}(\p(a))(=\{\GA\ap\in Y\st\p(a)\in\GA\ap\}).
\end{equation}
Indeed, let $\GA\ap\in f_\p\inv(g_{\unlb}(a))$. Then $f_\p(\GA\ap)\in g_{\unlb}(a)$. Thus $a\in
\p\inv(\GA\ap)$, i.e. $\p(a)\in\GA\ap$. So, $\GA\ap\in g_{\unlb\ap}(\p(a))$. Hence $f_\p\inv(g_{\unlb}(a))\sbe g_{\unlb\ap}(\p(a))$. Conversely,
let $\GA\ap\in g_{\unlb\ap}(\p(a))$, i.e. $\p(a)\in\GA\ap$. Then $a\in\p\inv(\GA\ap)=f_\p(\GA\ap)$. Hence $f_\p(\GA\ap)\in g_{\unlb}(a)$. Then $\GA\ap\in f_\p\inv(g_{\unlb}(a))$. So,
$f_\p\inv(g_{\unlb}(a))\spe g_{\unlb\ap}(\p(a))$. Thus the equation (\ref{contfphic}) is verified and
we get that $f_\p$ is a continuous function.

Let us now show that $f_\p(Y_0)\sbe X_0$. Let $u\ap\in Y_0$. Then $u\ap\in Ult(B)$. Hence $f_\p(u\ap)=\p\inv(u\ap)\in Ult(A)=X_0$. Therefore, $f_\p(Y_0)\sbe X_0$.

Let $x,y\in Y_0$ and $xR\ap y$. Then $x,y\in Ult(B\ap)$. Hence $f_\p(x),f_\p(y)\in Ult(B)$. We will show that $f_\p(x)R f_\p(y)$. We have that $f_\p(x)R f_\p(y)\iff f_\p(x)\times f_\p(y)\sbe C\iff(\fa a\in\p\inv(x))(\fa b\in\p\inv(y))(aCb)$. Now, since $xR\ap y$, we have that $x\times y\sbe C\ap$, i.e. $(\fa a\in x)(\fa b\in y)(aC\ap b)$. Thus, if $a\in\p\inv(x)$ and $b\in\p\inv(y)$, then $\p(a)\in x$ and $\p(b)\in y$; hence $\p(a)C\ap\p(b)$ and, therefore, $aCb$. This shows that $f_\p(x)R f_\p(y)$.
So,
$$G^a(\p)\in\PCS(G^a(\unlb\ap),G^a(\unlb)).$$

Clearly, for every precontact algebra $(B,C)$, $G^a(id_{(B,C)})=id_{G^a(B,C)}$. Let $\p\in\PCA((B,C),(B_1,C_1))$ and $\psi\in\PCA((B_1,C_1),(B_2,C_2))$. Let $f_\p=G^a(\p)$, $f_\psi=G^a(\psi)$ and $G^a(B_2,C_2)=(Z,Z_0,R)$. Then, for every $\GA\in Z$, we have that $G^a(\psi\circ\p)(\GA)=(\psi\circ\p)\inv(\GA)=\p\inv(\psi\inv(\GA))=f_\p(f_\psi(\GA))=(G^a(\p)\circ G^a(\psi))(\GA)$. We get
 that $G^a$ is a contravariant functor.

Let $\unlx=(X,X_0, R)\in\card{\PCS}$. We will fix a $\PCS$-isomorphism $t_{\unlx}$ between $\unlx$ and $G^a(G^t(\unlx)$. We have that $G^t(X,X_0,R)=\unlb(\unlx)=(RC(X,X_0),C_{\unlx}).$ Set
$$\unlb=(B,C)=(RC(X,X_0),C_{\unlx}).$$
Let $G^a(\unlb)=(Y,Y_0,R\ap)$. Then $Y=Clans(\unlb)$,  $Y_0=Ult(B)$ and, for every $x,y\in Y_0$, $xR\ap y\iff x\times y\sbe C$.
By \cite[Proposition 4.7(a)]{DV11} and Lemma \ref{isombool}, the contact algebras
\begin{equation}\label{contfphicnovo}
(CO(X_0),(C_R)^{\#})\ \mbox{ and }\ (RC(X,X_0),C_{(X,X_0)})\ \mbox{ are isomorphic.}
\end{equation}
Hence there exists a bijection between the sets
$$Clans(CO(X_0),(C_R)^{\#})\ \mbox{ and }\ Clans(RC(X,X_0),C_{(X,X_0)}).$$
 We have that for every $F,G\in CO(X_0)$,
$$F(C_R)^{\#}G\iff(FC_RG)\vee(GC_RF)\vee(F.G\neq 0)\iff$$
$$\iff[(\ex x\in F)(\ex y\in G)(xRy)]\vee[(\ex y\in G)(\ex x\in F)(yRx)]\vee(F.G\neq 0)\iff$$
$$\iff \cl_X(F)(C_{\unlx})^{\#}\cl_X(G).$$
 Hence the contact algebras
 \begin{equation}\label{contfphicnovo1}
 (RC(X,X_0),(C_{\unlx})^\#),\ \  (CO(X_0),(C_R)^{\#})\  \mbox{ and } \ (RC(X,X_0),C_{(X,X_0)})
 \end{equation}
  are isomorphic. In fact, the contact algebras
  \begin{equation}\label{contfphicnovo2}
(RC(X,X_0),(C_{\unlx})^\#)\ \mbox{ and }\ (RC(X,X_0),C_{(X,X_0)})\ \mbox{ coincide.}
  \end{equation}
 Therefore, the sets $$Clans(RC(X,X_0),C_{\unlx})(=Y)\  \mbox{ and }\ Clans(RC(X,X_0),C_{(X,X_0)})$$ coincide.
By \cite[Proposition 4.1(ii)]{DV1}, for every $x\in X$, $\s_x\in Clans(RC(X),C_X)$. Hence, by Lemma \ref{obvle}, for every subalgebra  $ A $  of the Boolean algebra $RC(X)$,
\begin{equation}\label{sxb}
 \s_x^A\in Clans(A,C_X\cap A^2).
 \end{equation}
 Therefore, for every $x\in X$,
 \begin{equation}\label{sxb1}
 \s_x^{RC(X,X_0)}\in Clans(RC(X,X_0),C_{(X,X_0)}).
 \end{equation}
So, the following map is well-defined:
\begin{equation}\label{txcnc}
t_{\unlx}:(X,X_0,R)\lra G^a(G^t(X,X_0,R))(=(Y,Y_0,R\ap)), \ \ x\mapsto \s_x^{RC(X,X_0)};
\end{equation}
we will  show that it is a homeomorphism.  We start by proving that $t_{\unlx}$ is a surjection.
 Let $\GA\in Y$. Then $\GA\in Clans(\unlb)$. From the above considerations, we get that $\GA\ap=r_{X_0,X}(\GA)\in Clans(CO(X_0),C_R)$. Hence, by (PCS5), there exists $x\in X$ such that $\GA\ap=\GA_{x,X_0}$. Since, by Lemma  \ref{isombool}, $\GA=e_{X_0,X}(\GA\ap)$, we get that $\GA=\s_x^B=t_{\unlx}(x)$. So, $t_{\unlx}$ is a surjection. For showing that $t_{\unlx}$ is a injection, let $x,y\in X$ and $x\neq y$. Since $X$ is a $T_0$-space, there exists an open subset $U$ of $X$ such that $|U\cap\{x,y\}|=1$. We can suppose,
 without loss of generality, that $x\in U$ and $y\nin U$. Since $B$ is a closed base of $X$, there exists $F\in B$ such that $x\in X\stm F\sbe U$. Then $y\in F$ and $x\nin F$. Hence $F\in \s_y^B$ and $F\nin\s_x^B$, i.e. $t_{\unlx}(x)\neq t_{\unlx}(y)$. So, $t_{\unlx}$ is an injection. Thus $t_{\unlx}$ is a bijection.
We will now prove that $t_{\unlx}$ is a continuous map. We have that the family $\{g_{\unlb}(F)=\{\GA\in Y\st F\in\GA\}\st F\in B\}$ is a closed base of $Y$. Let $F\in B$. We will show that
\begin{equation}\label{txcontc}
t_{\unlx}\inv(g_{\unlb}(F))=F.
\end{equation}
Let $x\in F$. Set $t_{\unlx}(x)=\GA$. Then $\GA=\s_x^B$. Since $F\in\GA$, we get that $\GA\in g_{\unlb}(F)$. Thus $t_{\unlx}(F)\sbe g_{\unlb}(F)$, i.e. $F\sbe t_{\unlx}\inv(g_{\unlb}(F))$. Conversely, let $x\in t_{\unlx}\inv(g_{\unlb}(F))$. Set $\GA=t_{\unlx}(x)$. Then $\GA\in g_{\unlb}(F)$. Hence $F\in\GA$. Since $\GA=\s_x^B$, we get that $x\in F$. Hence $F\spe t_{\unlx}\inv(g_{\unlb}(F))$. So, $F= t_{\unlx}\inv(g_{\unlb}(F))$. This shows that $t_{\unlx}$ is a continuous map. For showing that $t_{\unlx}\inv$ is a continuous map, let $F\in B$. Using (\ref{txcontc}) and the fact that $t_{\unlx}$ is a bijection, we get that $t_{\unlx}(F)=g_{\unlb}(F)$. Hence $(t_{\unlx}\inv)\inv(F)=g_{\unlb}(F)$. This shows that $t_{\unlx}\inv$ is a continuous map. So, $t_{\unlx}$ is a homeomorphism.

We will now  show that $t_{\unlx}(X_0)=Y_0$. Let $x\in X_0$. Set $\GA=t_{\unlx}(x)$. Then $\GA=\s_x^B$ and $r_{X_0,X}(\GA)=\{F\in CO(X_0)\st x\in F\}=u_x^{X_0}\in Ult(CO(X_0))$. Then, by Lemma \ref{isombool}, $\GA=e_{X_0,X}(u_x^{X_0})\in Ult(B)=Y_0$. Hence $t_{\unlx}(X_0)\sbe Y_0$. Let now $\GA\in Y_0$. Then $\GA\in Ult(B)$ and thus $u=r_{X_0,X}(\GA)\in Ult(CO(X_0))$. Clearly, there exist $x\in X_0$ such that $u=u_x^{X_0}$. Then $\GA=e_{X_0,X}(u_x^{X_0})=\s_x^B=t_{\unlx}(x)$. Therefore, $t_{\unlx}(X_0)\spe Y_0$. We have proved that $t_{\unlx}(X_0)= Y_0$.

Finally, we will prove that for every $x,y\in X_0$, $xRy\iff t_{\unlx}(x)R\ap t_{\unlx}(y)$. So, let $x,y\in X_0$ and $xRy$. Then $t_{\unlx}(x)=\s_x^B$ and $t_{\unlx}(y)=\s_y^B$. We have that $$t_{\unlx}(x)R\ap t_{\unlx}(y)\iff\s_x^B\times\s_y^B\sbe C_{\unlx}\iff [((x\in F\in B)\we(y\in G\in B))\Rightarrow (FC_{\unlx}G)].$$ This shows that, since $xRy$, $t_{\unlx}(x)R\ap t_{\unlx}(y)$ holds. Conversely, let $x,y\in X_0$ and $t_{\unlx}(x)R\ap t_{\unlx}(y)$. Suppose that $x(-R)y$. Then there exist $F,G\in CO(X_0)$ such that $x\in F$, $y\in G$ and $R\cap (F\times G)=\ems$. This implies that $\cl_X(F)(-C_{\unlx})\cl_X(G)$, a contradiction. Hence, $xRy$.

All this shows that
\begin{equation}\label{txpcsisom}
t_{\unlx}  \mbox{ is a } \PCS\mbox{-isomorphism}.
\end{equation}

Let $\unlb=(B,C)$ be a precontact algebra and let us set $(X,X_0,R)=G^a(\unlb)$. Then $G^t(X,X_0,R)=(RC(X,X_0),C_{\unlx})$ and, using Theorem \ref{mth}(a),
we get that the map $$g_{\unlb}:\unlb\lra G^t(G^a(\unlb))(=(RC(X,X_0),C_{\unlx})), \ \  a\mapsto g_{\unlb}(a)=\{\GA\in X\st a\in\GA\},$$
is a PCA-isomorphism.
This implies immediately that $g_{\unlb}$ is also a $\PCA$-isomorphism.

We will now show that
%
$$t:Id_{\PCS}\lra   G^a\circ  G^t,$$
%
 defined by
%
 $t(\unlx)=\tcxx, \ \ \fa \unlx\in\card\PCS,$
%
is a natural isomorphism.

Let $\unlx=(X,X_0,R)$, $\unly=(Y,Y_0,R\ap)$, $\unlx,\unly\in\card{\PCS}$,
 $f\in\PCS(\unlx,\unly)$ and $\fs= G^a(G^t(f))$. We will show
that
\begin{equation}\label{natoft}
\fs\circ \tcxx=\tcyy\circ f.
\end{equation}
 Set $\p_f=G^t(f)$, $A=RC(X,X_0)$, $B=RC(Y,Y_0)$, $\unlb=G^t(\unly)(=(B,C_{\unly}))$.
Let $x\in X$. Then
$$(\tcyy\circ f)(x)=\tcyy(f(x))=\s_{f(x)}^B=\{F\in B\st f(x)\in B\}.$$
Further,
$\fs(\tcxx(x))=\fs( \s_x^A)=\p_f\inv(\s_x^A)$. Set $\GA\ap= \p_f\inv(\s_x^A)$.
Then $\GA\ap=\{G\in B\st \p_f(G)\in\s_x^A\}=\{\cl_Y(G_0)\st G_0\in CO(Y_0), x\in\cl_X(X_0\cap f\inv(G_0))\}$.
Let $G_0\in CO(Y_0)$ and $\cl_Y(G_0)\in\GA\ap$. Then $x\in\cl_X(X_0\cap f\inv(G_0))$ and thus $f(x)\in f(\cl_X(X_0\cap f\inv(G_0))\sbe \cl_Y(f(X_0\cap f\inv(G_0)))\sbe \cl_Y(G_0)$. Therefore, $$\GA\ap\sbe\s_{f(x)}^B.$$ We have that $\GA\ap\in Clans(\unlb)=Clans(RC(Y,Y_0),C_{(Y,Y_0)})$. Hence $\GA\ap_r=r_{Y_0,Y}(\GA\ap)\in Clans(CO(Y_0),C_{R\ap})$. Thus, by (PCS5), there exists $y\in Y$ such that $\GA\ap_r=\GA_{y,Y_0}$. Then
\begin{equation}\label{gaap}
\GA\ap=\s_y^B.
\end{equation}
Since $\GA\ap\sbe\s_{f(x)}^B$ and $\GA\ap$ is a grill, we get, by Lemma \ref{pro41dv1}, that $\nu_{f(x)}^B\sbe\GA\ap$. According to \cite[Proposition 4.1]{DV1}, $\nu_{f(x)}^B$ is a filter of $B$. Hence, by Lemma \ref{grilllemma}, there exists an ultrafilter $u$ of $B$ such that $\nu_{f(x)}^B\sbe u\sbe \GA\ap$. Then $u\sbe \s_y^B$ and since $u$ is a grill of $B$, Lemma \ref{pro41dv1} implies that $\nu_y^B\sbe u$. So, we obtained that $\nu_{f(x)}^B\cup\nu_y^B\sbe u\sbe \GA\ap$. Then, for every $F\ap\in\nu_{f(x)}^B$ and every $G\ap\in\nu_y^B$, we have that $F\ap . G\ap\neq 0$, i.e. $\cl_Y(\int_Y(F\ap\cap G\ap))\nes$. Hence $\int_Y(F\ap\cap G\ap)\nes$ and thus $\int_Y(F\ap)\cap \int_Y(G\ap)\nes$, for every $F\ap\in\nu_{f(x)}^B$ and every $G\ap\in\nu_y^B$. Since $Y$ is a $T_0$-space, using \cite[Lemma 8.14]{DV11}, we get that $y=f(x)$. Therefore
$\GA\ap=\s_{f(x)}^B$. Thus $\fs\circ \tcxx=\tcyy\circ f$ and hence $t$ is a natural isomorphism.

Finally, we will prove that
$$g: Id_{\PCA}\lra  G^t\circ  G^a,\mbox{ where } g(\unlb)=g_{\unlb}, \ \ \fa \unlb\in\card\PCA,$$
%
 is a natural
isomorphism.

Let $\unla=(A,C)$, $\unlb=(B,C\ap)$, $\unla,\unlb\in\card{\PCA}$, $\p\in\PCA(\unla,\unlb)$ and $\ps=G^t(G^a(\p))$. We have
to prove that $g_{\unlb}\circ\p=\ps\circ g_{\unla}$. Set
$f=G^a(\p)$, $(X,X_0,R)=G^a(\unla)$ and $(Y,Y_0,R\ap)=G^a(\unlb)$. Then $\ps=G^t(f)(=\p_f)$. Let
$a\in A$. Then $g_{\unlb}(\p(a))=\{\GA\ap\in Y\st \p(a)\in\GA\ap\}.$
Further,
using Lemma \ref{isombool}, we get that
$$g_{\unlb}(\p(a))=\cl_Y(s_B(\p(a)))\ \ \mbox{ and }\ \  g_A(a)=\cl_X(s_A(a))$$
(see (\ref{sofa}) for the notation $s_A$).
Thus $$\ps(g_{\unla}(a))=\cl_Y(Y_0\cap f\inv(s_A(a))).$$
Let $u\ap\in Y_0\cap f\inv(s_A(a))$. Then $u\ap\in Ult(B)$ and $f(u\ap)\in s_A(a)$. Hence $\p\inv(u\ap)\in s_A(a)=\{u\in Ult(A)\st a\in u\}$. Thus $a\in\p\inv(u\ap)$, i.e. $\p(a)\in u\ap$. Therefore $u\ap\in s_B(\p(a))$. So, $Y_0\cap f\inv(s_A(a))\sbe s_B(\p(a))$. Conversely, let $u\ap\in s_B(\p(a))$. Then $u\ap\in Y_0$ and $\p(a)\in u\ap$. Hence $a\in\p\inv(u\ap)=f(u\ap)$. Thus $f(u\ap)\in s_A(a)$. Therefore, $u\ap\in Y_0\cap f\inv(s_A(a))$. So, $Y_0\cap f\inv(s_A(a))\spe s_B(\p(a))$ and we get that $Y_0\cap f\inv(s_A(a))= s_B(\p(a))$. Hence $\ps(g_{\unla}(a))=\cl_Y(s_B(\p(a)))=g_{\unlb}(\p(a))$.
 So,
$g$ is a natural isomorphism.

We have proved that $(G^t,G^a,g,t)$ is a duality between the categories $\PCS$ and $\PCA$.
 \sqs

 In the proof of the next corollary  (as well as in the proofs of all further corollaries of Theorem \ref{maindualitytheorem}), we will use the notation from (the proof of) Theorem \ref{maindualitytheorem}.

\begin{cor}\label{corgtnew}
Let $\unlx=(X,X_0,R)$, $\unly=(Y,Y_0,R\ap)$ and $f\in\PCS(\unlx,\unly)$. Then, for every $H\in RC(Y,Y_0)$,
\begin{equation}\label{gtnewdef}
G^t(f)(H)=f\inv(H)
\end{equation}
 (see (\ref{deftetafc}) for the notation $G^t$).
\end{cor}

\doc  Recall that $RC(Y,Y_0)=\{\cl_Y(F)\st F\in CO(Y_0)\}$. Let now $H\in RC(Y,Y_0)$ and $\fs= G^a(G^t(f))$. Then, by (\ref{natoft}), $\fs\circ \tcxx=\tcyy\circ f$. Thus,
$$f\inv(H)=(\tcxx)\inv(\fs\inv(\tcyy(H))).$$
Let $\unlb=G^t(\unlx)$, $\unlb\ap=G^t(\unly)$ and $\p=G^t(f)$.
Using (\ref{txpcsisom}), (\ref{txcontc}) and (\ref{contfphic}), we get that $\tcyy(H)=g_{\unlb\ap}(H)$ (see (\ref{fg}) for the notation $g_{\unlb\ap}(H)$), $\fs\inv(\tcyy(H))=\fs\inv(g_{\unlb\ap}(H))=g_{\unlb}(\p(H))$ and $(\tcxx)\inv(\fs\inv(\tcyy(H)))=(\tcxx)\inv(g_{\unlb}(\p(H)))=\p(H)$. Hence, $f\inv(H)=G^t(f)(H)$.
\sqs

As it was already noted, Theorem \ref{maindualitytheorem} implies the Stone Duality Theorem \cite{ST,kop89} and its connected version \cite{DV11}.

\begin{cor}\label{corstonedth}{\rm (\cite{ST,kop89})}
The category  $\Bool$  and the category  $\Stone$ of compact zero-dimensional Hausdorff spaces and continuous maps
are dually equivalent.
\end{cor}

\doc  Let $\B$ be the full subcategory of the category $\PCA$  having as objects all (pre)contact algebras of the form $(B,\rho_s^B)$ (see Example \ref{extrcr} for the notation $\rho_s^B$). It is easy to see that if $A$ and $B$ are Boolean algebras, then a function $\p:A\lra B$ induces a $\PCA$-morphism between $(A,\rho_s^A)$ and $(B,\rho_s^B)$ iff $\p$ is a Boolean homomorphism. Hence the categories $\Bool$ and $\B$ are isomorphic.

Let $\S$ be the full subcategory of the category $\PCS$ having as objects all 2-precontact spaces of the form $(X,X,D_X)$, where $D_X$ is the diagonal of $X$. Obviously, if $(X,X,D_X)\in\card{\S}$ then $X$ is a Stone space. Conversely, it is easy to see that every triple of the form $(X,X,D_X)$, where $X$ is a Stone space, is a 2-precontact space.
 Thus the objects of the category $\S$ are all triples of the form $(X,X,D_X)$, where $X$ is a Stone space. Clearly, if $X,Y\in\card{\Stone}$, then a function $f:X\lra Y$ induces a $\PCS$-morphism between $(X,X,D_X)$ and $(Y,Y,D_Y)$ iff $f$ is a continuous map. Hence the categories $\S$ and $\Stone$ are isomorphic.

 We will now show that $G^a(\B)\sbe\S$ and $G^t(\S)\sbe\B$. Indeed, if $B$ is a Boolean algebra, then $G^a(B,\rho_s^B)=(Clans(B,\rho_s^B), Ult(B),R)$, where, for every $u,v\in Ult(B)$, $uRv\iff u\times v\sbe\rho_s^B$. Obviously, we have that for every $u,v\in Ult(B)$, $uRv\iff u=v$. Further, by \cite[Example 3.1]{DV1}, $Clans(B,\rho_s^B)=Ult(B)$. Also, by Theorem \ref{maindualitytheorem}, $G^a(B,\rho_s^B)$ is a 2-precontact space. Hence $G^a(B,\rho_s^B)=(X,X,D_X)$, where $X\in\card{\Stone}$.  So, $G^a(\B)\sbe\S$.  Let now $X\in\card{\Stone}$ and $\unlx=(X,X,D_X)$. Then $G^t(\unlx)=(RC(X,X),C_{\unlx})$. We have that $RC(X,X)=CO(X)$ and for every $F,G\in CO(X)$, $FC_{\unlx}G\iff[(\ex x\in F)(\ex y\in G)(xD_X y)]\iff F\cap G\nes\iff F.G\neq 0\iff F\rho_s^{CO(X)} G$. Thus $G^t(\unlx)\in\card{\S}$. So, we get that $G^t(\S)\sbe\B$.
Therefore we obtain, using Theorem \ref{maindualitytheorem}, that the restriction of the contravariant functor $G^a$ to the category $\B$ is a duality between the categories $\B$ and $\S$. This implies that the categories $\Bool$ and $\Stone$ are dually equivalent.
\sqs

Now we will obtain as a corollary of Theorem \ref{maindualitytheorem} the connected version of the Stone Duality Theorem proved in \cite{DV11}.

\begin{cor}\label{constonedualitycor}{\rm (\cite{DV11})}
The categories $\Boo$ and $\ECS$ are dually equivalent.
\end{cor}

\doc  Let $\B\ap$ be the full subcategory of the category $\PCA$  having as objects all (pre)contact algebras of the form $(B,\rho_l^B)$ (see Example \ref{extrcr} for the notation $\rho_l^B$). It is easy to see that if $A$ and $B$ are Boolean algebras, then a function $\p:A\lra B$ induces a $\PCA$-morphism between $(A,\rho_l^A)$ and $(B,\rho_l^B)$ iff $\p$ is a Boolean homomorphism. Hence the categories $\Bool$ and $\B\ap$ are isomorphic.

Let $\S\ap$ be the full subcategory of the category $\PCS$ having as objects all 2-precontact spaces of the form $(X,X_0,(X_0)^2)$. We will show that if $(X,X_0,(X_0)^2)\in\card{\S\ap}$ then $(X,X_0)$ is a Stone 2-space. Indeed, set $B=CO(X_0)$ and $R=(X_0)^2$. Then, for every $F,G\in B$, we have that $FC_RG\iff [(\ex x\in F)(\ex y\in G)((x,y)\in R)]\iff [(F\nes)\we(G\nes)]\iff F\rho_l^{B}G$. Thus, using \cite[Example 3.1]{DV1}, we get that $Clans(B,C_R)=Clans(B,\rho_l^B)=Grills(B)$. Now it becomes obvious that the axiom (S2S4) is fulfilled because the axiom (PCS5) is fulfilled. The rest is clear, so that $(X,X_0)$ is a Stone 2-space. Arguing similarly, we get that, conversely, if $(X,X_0)$ is a Stone 2-space then $(X,X_0,(X_0)^2)$ is a 2-precontact space. Hence, the objects of the category $\S\ap$ are
all triples $(X,X_0,(X_0)^2)$, where $(X,X_0)$ is a Stone 2-space (i.e., $(X,X_0)\in\card{\ECS}$).
 Clearly, if $(X,X_0),(Y,Y_0)\in\card{\ECS}$, then a function $f:X\lra Y$ induces a $\PCS$-morphism between $(X,X_0,X_0^2)$ and $(Y,Y_0,Y_0^2)$ $\iff$ ($f$ is a continuous map and $f(X_0)\sbe Y_0$) $\iff$ ($f\in\ECS((X,X_0),(Y,Y_0))$). Hence the categories $\S\ap$ and $\ECS$ are isomorphic.

 We will now show that $G^a(\B\ap)\sbe\S\ap$ and $G^t(\S\ap)\sbe\B\ap$.
 Indeed, if $B$ is a Boolean algebra, then $G^a(B,\rho_l^B)=(Clans(B,\rho_l^B), Ult(B),R)$, where, for every $u,v\in Ult(B)$, $uRv\iff u\times v\sbe\rho_l^B$. Since the ultrafilters do not contain the zero element, we get that $uRv$ for every $u,v\in Ult(B)$, i.e. $R=(Ult(B))^2$. Also, by Theorem \ref{maindualitytheorem}, $G^a(B,\rho_l^B)$ is a 2-precontact space. Hence $G^a(B,\rho_l^B)\in\card{\S\ap}$.   So, $G^a(\B\ap)\sbe\S\ap$.  Let now $(X,X_0)\in\card{\ECS}$ and $\unlx=(X,X_0,X_0^2)$. Then $G^t(\unlx)=(RC(X,X_0),C_{\unlx})$. We have that  for every $F,G\in CO(X_0)$, $(\cl_X(F)C_{\unlx}\cl_X(G))\iff[(\ex x\in F)(\ex y\in G)((x,y)\in X_0^2)]\iff [(F\nes)\we(G\nes)]\iff (\cl_X(F)\rho_l^{RC(X,X_0)} \cl_X(G))$. Thus $G^t(\unlx)\in\card{\S\ap}$. So, we get that $G^t(\S\ap)\sbe\B\ap$.
Therefore we obtain, using Theorem \ref{maindualitytheorem}, that the restriction of the contravariant functor $G^a$ to the category $\B\ap$ is a duality between the categories $\B\ap$ and $\S\ap$. This implies that the categories $\Bool$ and $\ECS$ are dually equivalent.
\sqs

Recall that: (a)(\cite{DV1}) if $(B,C)$ is a contact algebra, then a non-empty set $U$ of ultrafilters of $B$ is called a {\em clique in} $(B,C)$  if for any two $u,v\in U$ we have that $u\times v\sbe C$; (b)(\cite{DV1}) a subset of $B$ is a clan in $(B,C)$ iff it is a union of the elements of a clique in $(B,C)$.

\begin{cor}\label{faithfullm}
Let $\unlx=(X,X_0,R_1)$, $\unly=(Y,Y_0,R_2)$ and $\unlx,\unly\in |\PCS|$. If $f,g\in\PCS(\unlx,\unly)$ and $f|_{X_0}=g|_{X_0}$ then $f=g$.
\end{cor}

\doc  Let $B=RC(X,X_0)$, $B\ap=RC(Y,Y_0)$ and $$G^a(G^t(\unlx))=(X\ap,X_0\ap,R_1\ap), \ \ G^a(G^t(\unly))=(Y\ap,Y_0\ap,R_2\ap).$$ Let $f\ap=G^a(G^t(f))$. Then, by Theorem \ref{maindualitytheorem}, we have that $t_{\unly}\circ f=f\ap\circ t_{\unlx}$, i.e. $$f\ap(\s_x^B)=\s_{f(x)}^{B\ap}.$$

Let $x\in X\stm X_0$. Then $t_{\unlx}(x)=\s_x^B$ is a clan in $(B,(C_{\unlx})^{\#})$ (see (\ref{sxb1}) and (\ref{contfphicnovo2})). Hence it is a union of the elements of a clique $\{\s_{x_\a}^B\st x_\a\in X_0, \a\in A\}$ in $(B,(C_{\unlx})^{\#})$. Let $\GA\ap=\{\cl_Y(G\ap)\st G\ap\in CO(Y_0), f(x_\a)\in G\ap$ for some  $\a\in A\}$.  We will show that
\begin{equation}\label{imeq}
\s_{f(x)}^{B\ap}=\GA\ap.
\end{equation}
Indeed, let $G\ap\in CO(Y_0)$ and $f(x_\a)\in G\ap$ for some $\a\in A$. Set $G=(f|_{X_0})\inv(G\ap)$. Then $x_\a\in G\in CO(X_0)$. Thus $x\in\cl_X(G)$. Since $f$ is a continuous function, we get that $f(x)\in\cl_Y(f(G))\sbe\cl_Y(G\ap)$. Thus $\cl_Y(G\ap)\in\s_{f(x)}^{B\ap}$. Therefore $\GA\ap\sbe \s_{f(x)}^{B\ap}$. Set, for every $\a\in A$, $u_\a=t_{\unlx}(x_\a)(=\s_{x_\a}^B)$ and $v_\a=t_{\unly}(f(x_\a))(=\s_{f(x_\a)}^{B\ap})$. Obviously, each $u_\a$ is an ultrafilter of $B$. Since $f(X_0)\sbe Y_0$, we obtain that each $v_\a$ is an ultrafilter of $B\ap$. Also, $u_\a\times u_\b\sbe (C_{\unlx})^{\#}$ for all $\a,\b\in A$, i.e. $u_\a R_1\ap u_\b$ for all $\a,\b\in A$. Since $t_{\unlx}$ is a $\PCS$-isomorphism, we get that $x_\a R_1 x_\b$ for all $\a,\b\in A$. Then $f(x_\a) R_2 f(x_\b)$ for all $\a,\b\in A$, and thus $t_{\unly}(f(x_\a)) R_2\ap t_{\unly}(f(x_\b))$ for all $\a,\b\in A$. This means that $v_a\times v_\b\sbe (C_{\unly})^{\#}$ for all $\a,\b\in A$, and, hence, $\{v_\a\st \a\in A\}$ is a clique in  $(B\ap,(C_{\unly})^{\#})$. Thus $\GA\ap$ is a clan in $(B\ap,(C_{\unly})^{\#})$, i.e. $\GA\ap\in Y\ap$. Then there exists $y\in Y$ such that $\GA\ap=t_{\unly}(y)$, i.e. $\GA\ap=\s_y^{B\ap}$. Now we show that $y=f(x)$ exactly as  in the paragraph immediately after (\ref{gaap}). So, (\ref{imeq}) is proved.
Analogously we get that $\s_{g(x)}^{B\ap}=\{\cl_Y(G\ap)\st G\ap\in CO(Y_0), g(x_\a)\in G\ap\ \mbox{ for some }\ \a\in A\}$. Let $g\ap=G^a(G^t(g))$. Then Theorem \ref{maindualitytheorem} implies that $g\ap(\s_x^B)=\s_{g(x)}^{B\ap}.$ Since $f|_{X_0}=g|_{X_0}$, we get that $f\ap(\s_x^B)=g\ap(\s_x^B)$. Then $f(x)=g(x)$. So, $f=g$.
\sqs

\begin{cor}\label{maindualitytheoremcor1}
 The category\/ $\PCS$ is equivalent to the category\/ $\SAS$ of all Stone adjacency spaces and all continuous maps $f:(X_0,R)\lra(X_0\ap,R\ap)$ between them such that, for every $x,y\in X_0$, $xRy$ implies $f(x)R\ap f(y)$.
\end{cor}

\doc Let $F^t:\PCS\lra\SAS$ be the functor defined by $$F^t(X,X_0,R)=(X_0,R)$$ on the objects of the category $\PCS$,  and by
$$F^t(f)=f|_{X_0}$$ for every $f\in\PCS((X,X_0,R),(Y,Y_0,R\ap))$. Clearly, $F^t$ is a functor. We will show that it is full, faithful and isomorphism-dense, i.e. that it is an equivalence functor.

Let $(X_0,R)\in |\SAS|$. Set $\unlb=(CO(X_0),C_R)$ and $\unlx=G^a(\unlb)$. According to Theorem \ref{maindualitytheorem}, $\unlx\in |\PCS|$. Let $\unlx=(X\ap,X_0\ap,R\ap)$ and $B=CO(X_0)$. Then, by the Stone duality theorem, $f:X_0\lra X_0\ap,\ x\mapsto u_x,$ is a homeomorphism  (recall that $u_x=\{P\in B\st x\in P\}$). We will show that for every $x,y\in X_0$, $xRy\iff f(x)R\ap f(y)$. Indeed, let $x,y\in X_0$ and $xRy$. We have that $f(x)R\ap f(y)\iff u_x R\ap u_y\iff u_x\times u_y\sbe C_R\iff[(\fa P\in u_x)(\fa Q\in u_y)(PC_R Q)]\iff([(x\in P\in B)\we(y\in Q\in B)]\Rightarrow(R\cap (P\times Q)\nes))$. Now, since $xRy$, we obtain that $f(x)R\ap f(y)$. Conversely, let $x,y\in X_0$ and $f(x)R\ap f(y)$. Suppose that $x(-R)y$. Since $R$ is a closed relation, there exist $P,Q\in B$ such that $x\in P$, $y\in Q$ and $R\cap (P\times Q)=\ems$. We get that $f(x)(-R\ap) f(y)$, a contradiction. Hence, $xRy$. Therefore, $f$ is a $\SAS$-isomorphism between $(X_0,R)$ and $(X_0\ap,R\ap)$. Since $F^t(\unlx)=F^t(X\ap,X_0\ap,R\ap)=(X_0\ap,R\ap)$, we obtain that $F^t$ is isomorphism-dense.

Clearly,  Corollary \ref{faithfullm} implies that $F^t$ is a faithful functor.

Let $\unlx=(X,X_0,R_0)$, $\unly=(Y,Y_0,R_1)$ and $\unlx,\unly\in |\PCS|$. We will show that the restriction of $F^t$ to the hom-set $\PCS(\unlx,\unly)$ is a surjection (i.e., the correspondence $\PCS(\unlx,\unly)\lra \SAS(F^t(\unlx),F^t(\unly)), \ \  f\mapsto f|_{X_0},$ is a surjection). Let $g\in\SAS((X_0,R_0),(Y_0,R_1))$ and $\p=S^t(g):CO(Y_0)\lra CO(X_0),\ \ G\mapsto g\inv(G)$, where $S^t:\Stone\lra\Bool$ is the Stone duality functor. In fact, $$\p\in\PCA((CO(Y_0),C_{R_1}),(CO(X_0),C_{R_0})).$$ Indeed, let $G_1,G_2\in CO(Y_0)$ and $\p(G_1)C_{R_0}\p(G_2)$. Then there exist $x\in\p(G_1)$ and $y\in\p(G_2)$ such that $xR_0 y$. We get that $x\in g\inv(G_1)$, $y\in g\inv(G_2)$ and  $xR_0 y$. Thus $g(x)\in G_1$, $g(y)\in g(G_2)$ and $g(x)R_1 g(y)$. Hence, $G_1 C_{R_1} G_2$. So, $\p$ is a $\PCA$-morphism. Then Theorem \ref{maindualitytheorem} implies that $$G^a(\p)\in\PCS(G^a(CO(X_0),C_{R_0}),G^a(CO(Y_0),C_{R_1})).$$ Let $G^a(CO(X_0),C_{R_0})=(X\ap, X_0\ap,R_0\ap)$ and $G^a((CO(Y_0),C_{R_1})=(Y\ap,Y_0\ap,R_1\ap)$. Then $X_0\ap=S(S^t(X_0))$ and $Y_0\ap=S(S^t(Y_0))$, where $S:\Bool\lra \Stone$ is the Stone duality functor. Thus $f_{X_0}:X_0\lra X_0\ap, \ x\mapsto u_x,$ and $f_{Y_0}:Y_0\lra Y_0\ap, \ y\mapsto u_y,$ are homeomorphisms and $f_{Y_0}\circ g=S(S^t(g))\circ f_{X_0}$. Let $g\ap=S(S^t(g))$. Then we get that $u_{g(x)}=g\ap(u_x)$, for every $x\in X_0$. We will show that $g\ap=(G^a(\p))|_{X_0\ap}$. Indeed, for every $x\in X_0$, we have that $(G^a(\p))(u_x)=\p\inv(u_x)=\{G\in CO(Y_0)\st \p(G)\in u_x\}=\{G\in CO(Y_0)\st x\in\p(G)\}=\{G\in CO(Y_0)\st x\in g\inv(G)\}=\{G\in CO(Y_0)\st g(x)\in G\}=u_{g(x)}=g\ap(u_x)$. Hence, $F^t(G^a(\p))=g\ap$. Since $(CO(X_0),C_{R_0})$ and $(RC(X,X_0),C_{\unlx})$ are $\PCA$-isomorphic and, analogously, $(CO(Y_0),C_{R_1})$ and $(RC(Y,Y_0),C_{\unly})$ are $\PCA$-isomorphic, we get that $G^a(CO(X_0),C_{R_0})$ is $\PCS$-isomorphic to $G^a(G^t(\unlx))$ and, analogously, $G^a(CO(Y_0),C_{R_1})$ is $\PCS$-isomorphic to $G^a(G^t(\unly))$. Now, using Theorem \ref{maindualitytheorem}, we obtain that there exists $g_1\in\PCS(\unlx,\unly)$ such that $F^t(g_1)=g$. Thus, the functor $F^t$ is full.

All this shows that $F^t:\PCS\lra\SAS$ is an equivalence functor.
\sqs

\begin{cor}\label{maindualitytheoremcor123}{\rm (\cite{BBSV})}
The categories $\PCA$ and $\SAS$ are dually equivalent.
\end{cor}

\doc Clearly, it follows from Theorem \ref{maindualitytheorem} and Corollary \ref{maindualitytheoremcor1}.
\sqs

\begin{rem}\label{bbsvrem}
\rm
It is not difficult to see that a direct proof of Corollary \ref{maindualitytheoremcor123} can be easily obtained  using Theorem \ref{th3}.

Note that Corollary \ref{maindualitytheoremcor123} is formulated in \cite{BBSV} using the relation {\em non-tangential inclusion}\/ (defined here in  (\ref{nontangincl}); see also   Remark \ref{ctrdiez}) instead of the precontact relation; the relation non-tangential inclusion is called there {\em subordination}.
\end{rem}

\begin{cor}\label{maindualitytheoremcor111}
For every Stone adjacency space $(X_0,R)$ there exists a unique (up to homeomorphism) topological space $X$ such that the triple $(X,X_0,R)$ is a 2-precontact space (and, thus, $X$ is a compact semiregular $T_0$-space).
\end{cor}

\doc By Corollary \ref{maindualitytheoremcor1}, there exists an equivalence functor $F^s:\SAS\lra\PCS$ such that the compositions $F^t\circ F^s$ and $F^s\circ F^t$ are naturally isomorphic to the corresponding identity functors. Let $(X_0,R)\in|\SAS|$. Since $F^t$ is isomorphism-dense, there exists $(X\ap,X_0\ap, R\ap)\in|\PCS|$ such that $F^t(X\ap,X_0\ap, R\ap)$ is $\SAS$-isomorphic to $(X_0,R)$. Then, clearly, there exists $\unlx=(X,X_0,R)\in|\PCS|$. Suppose that there exists $\unlx_1=(X_1,X_0,R)\in|\PCS|$. Then $F^s(F^t(\unlx))=F^s(F^t(\unlx_1))$ and thus $\unlx$ is $\PCS$-isomorphic to $\unlx_1$. This implies that $X$ is homeomorphic to $X_1$. Note that, by Definition \ref{precontactspace} and Corollary \ref{cormth},   $X$ is a compact semiregular $T_0$-space.
\sqs

\begin{cor}\label{maindualitytheoremcor2}{\sc (A Duality Theorem for Contact Algebras)}
The full subcategory\/ $\CA$ of the category\/ $\PCA$ whose objects are all contact algebras is dually equivalent to the category\/  $\CS$ of all 2-contact spaces and all continuous maps $f:(X,X_0)\lra(X\ap,X_0\ap)$ between them such that $f(X_0)\sbe X_0\ap$.
\end{cor}

\doc Let $\S''$ be the full subcategory of the category $\PCS$ whose objects are all 2-precontact spaces $(X,X_0,R)$ for which $R$ is a reflexive and symmetric relation.
Then the categories $S''$ and $\CS$ are isomorphic. Indeed, let $F^c:\S''\lra\CS$ be defined by $F^c(X,X_0,R)=(X,X_0)$ on the objects of the category $\S''$, and by $F^c(f)=f$ on the morphisms of the the category $\S''$. Then, using \cite[Proposition 7.7]{DV11}, we get that $F^c$ is well-defined and, obviously, it is a functor. By \cite[Lemma 7.5]{DV11}, for every $(X,X_0)\in|\CS|$ there exists a unique relation $R_{(X,X_0)}$ on $X_0$ such that $(X,X_0,R)\in|\S''|$; as it is shown in \cite[Lemma 7.5]{DV11}, the relation $R_{(X,X_0)}$ is defined as follows: for every $x,y\in X_0$,
\begin{equation}\label{xrycont}
xR_{(X,X_0)}y\iff ((\fa F\in u_x)(\fa G\in u_y)(\cl_X(F)\cap\cl_X(G)\nes)),
\end{equation}
where  $u_x=\{A\in CO(X_0)\st x\in A\}$ and analogously for $u_y$ (see (\ref{ux})).
Now we set $F^d(X,X_0)=(X,X_0,R_{(X,X_0)})$, for every $(X,X_0)\in|\CS|$. Also, for every $f\in\CS((X,X_0),(Y,Y_0))$, we set $F^d(f)=f$. Then $f\in\S''(F^d(X,X_0),F^d(Y,Y_0))$. Indeed, set $R=R_{(X,X_0)}$ and $R\ap=R_{(Y,Y_0)}$. Let $x,y\in X_0$ and $xRy$. Then $f(x),f(y)\in Y_0$ (since $f(X_0)\sbe Y_0$). Suppose that $f(x)(-R\ap) f(y)$. Then there exist $F\ap\in u_{f(x)}$ and $G\ap\in u_{f(y)}$ such that $\cl_Y(F\ap)\cap\cl_Y(G\ap)=\ems$. Let $F=X_0\cap f\inv(F\ap)$ and
$G=X_0\cap f\inv(G\ap)$. Then $F\in u_x$ and $G\in u_y$. Thus there exists $z\in \cl_X(F)\cap\cl_X(G)$. Since $f$ is a continuous function, we get that $f(z)\in f(\cl_X(F))\sbe\cl_Y(f(F))\sbe\cl_Y(F\ap)$ and, analogously, $f(z)\in\cl_Y(G\ap)$. Therefore, $\cl_Y(F\ap)\cap\cl_Y(G\ap)\nes$, a contradiction. Hence $f(x) R\ap f(y)$.
So, $f\in\S''(F^d(X,X_0),F^d(Y,Y_0))$. All this shows that we have defined a functor $F^d:\CS\lra\S''$. Clearly, $F^c\circ F^d=Id_{\CS}$ and $F^d\circ F^c=Id_{{\S}''}$.
Hence the categories $\S''$ and $\CS$ are isomorphic.

We will now show that $G^a(\CA)\sbe\S''$ and $G^t(\S'')\sbe\CA$. The first inclusion follows from \cite[Lemma 3.5(d,e)]{DV11}, and the second one follows from Proposition \ref{prop1}(b) since, for every $\unlx=(X,X_0,R)\in |\PCS|$, the precontact algebras $(RC(X,X_0),C_{\unlx})$ and $(CO(X_0),C_R)$ are isomorphic. Now, applying Theorem \ref{maindualitytheorem}, we obtain that the categories $\CA$ and $\CS$ are dually equivalent.
\sqs

\begin{cor}\label{corgtnew1}
Let  $f\in\CS((X,X_0),(Y,Y_0))$. Then, for every $H\in RC(Y,Y_0)$,
\begin{equation}\label{gtnewdef1}
f\inv(H)=\cl_X(X_0\cap f\inv(H)).
\end{equation}
\end{cor}

\doc It follows from Corollary \ref{maindualitytheoremcor2}, Corollary \ref{corgtnew} and the obvious fact that  that $X_0\cap f\inv(H)=X_0\cap f\inv(H\cap Y_0)$ (because $f(X_0)\sbe Y_0$).
\sqs

The next corollary follows immediately from Corollary \ref{maindualitytheoremcor123} and Proposition \ref{prop1}:

\begin{cor}\label{maindualitytheoremcor123c}{\rm (\cite{BBSV})}
The category $\CA$ is dually equivalent to the full subcategory $\CSAS$ of the category $\SAS$ whose objects are all Stone adjacency spaces $(X,R)$ such that $R$ is a reflexive and symmetric relation.
\end{cor}

Let us derive it from Corollary \ref{maindualitytheoremcor2} as well. Indeed, by the proof of Corollary \ref{maindualitytheoremcor2}, the categories $\S''$ and $\CA$
are dually equivalent; by Corollary \ref{maindualitytheoremcor1}, the functor $F^t:\PCS\lra\SAS$ is an equivalence; obviously, $F^t(\S'')=\CSAS$ and hence the categories $\S''$ and $\CSAS$
are  equivalent; therefore, the categories $\CSAS$ and $\CA$
are dually equivalent.

\begin{cor}\label{maindualitytheoremcor3}{\sc (A Duality Theorem for Complete Contact Algebras)}
The full subcategory \/ $\CCA$ of the category\/ $\PCA$ whose objects are all complete contact algebras is dually equivalent to the category\/ $\CSRS$ of all C-semiregular spaces and all continuous maps  between them which preserve u-points.
\end{cor}

\doc Let $\S'''$ be the full subcategory of the category $\PCS$ whose objects are all 2-precontact spaces $(X,X_0,R)$ for which $R$ is a reflexive and symmetric relation and $X_0$ is extremally disconnected.
Then the categories $S'''$ and $\CSRS$ are isomorphic. Indeed, let $F^e:\S'''\lra\CSRS$ be defined by $F^e(X,X_0,R)=X$ on the objects of the category $\S'''$, and by $F^e(f)=f$ on the morphisms of the the category $\S'''$. Then, by \cite[Proposition 7.7]{DV11} and Lemma \ref{csem2cont}, $F^e$ is well-defined on the objects of the category $S'''$. Let $\unlx=(X,X_0,R)\in |\S'''|$, $\unly=(Y,Y_0,R\ap)\in |\S'''|$ and $f\in\S'''(\unlx,\unly)$. We will show that $F^e(f)\in\CSRS(F^e(\unlx),F^e(\unly))$. We have, by \cite[Proposition 7.7]{DV11}, Lemma \ref{csem2cont} and Theorem \ref{novocor2}, that $X_0=\{x\in X\st x \mbox{ is an u-point of }X\}$ and $Y_0=\{y\in Y\st y \mbox{ is an u-point of }Y\}$. Since $f$ is a continuous map and $f(X_0)\sbe Y_0$, we get that $f$ is a $\CSRS$-morphism. So, $F^e$ is well-defined.
Obviously, $F^e$ is a functor.
Let now $X\in|\CSRS|$.
 We set $F^f(X)=(X,X_0,R_{(X,X_0)})$, where $X_0=\{x\in X\st x \mbox{ is an u-point of }X\}$ and the relation $R_{(X,X_0)}$ on $X_0$ is defined by (\ref{xrycont}). Then,
 using Theorem \ref{novocor2} and the proof of Corollary \ref{maindualitytheoremcor2}, we get that $F^f(X)\in |\S'''|$.
 Further, for every $f\in\CSRS(X,Y)$, we set $F^f(f)=f$. Let $F^f(X)=(X,X_0,R_{(X,X_0)})$ and $F^f(Y)=(Y,Y_0,R_{(Y,Y_0)})$. Then, clearly, $f(X_0)\sbe Y_0$ and since, by
 Theorem \ref{novocor2},  $(X,X_0)$ and $(Y,Y_0)$ are 2-contact spaces, we obtain that $f\in\CS((X,X_0),(Y,Y_0))$; now, by  Corollary \ref{maindualitytheoremcor2}, $f=F^d(f)$ and we get that $f\in\S''(F^d(X,X_0),F^d(Y,Y_0))$; this implies that
  $f\in\S'''(F^f(X),F^f(Y))$. All this shows that we have defined a functor $F^f:\CSRS\lra\S'''$. Clearly, $F^e\circ F^f=Id_{\CSRS}$ and $F^f\circ F^e=Id_{{\S}'''}$.
Hence the categories $\S'''$ and $\CSRS$ are isomorphic.

We will now show that $G^a(\CCA)\sbe\S'''$ and $G^t(\S''')\sbe\CCA$. By  Corollary \ref{maindualitytheoremcor2}, we have that $G^a(\CCA)\sbe\S''$ and then Theorem \ref{mth}(a) together with Lemma \ref{lmrcx0}
 imply that even $G^a(\CCA)\sbe\S'''$. Further, by Corollary \ref{maindualitytheoremcor2}, $G^t(\S'')\sbe\CA$. Since, for every $\unlx=(X,X_0,R)\in |\S'''|$, we have that $CO(X_0)=RC(X_0)$ and thus $RC(X,X_0)=RC(X)$, we obtain that $G^t(\S''')\sbe\CCA$.
 Now, applying Theorem \ref{maindualitytheorem}, we obtain that the categories $\CCA$ and $\CSRS$ are dually equivalent.
\sqs

We are now going  to present Corollary \ref{maindualitytheoremcor2} in a form similar to that of Corollary \ref{maindualitytheoremcor3}. For doing this we need to recall two definitions from \cite{GG}. We first introduce a new notion.

\begin{defi}\label{mepair}
\rm
Let $X$ be a topological space and $B$ be a Boolean subalgebra of the Boolean algebra $RC(X)$. Then the pair $(X,B)$ is called a {\em mereotopological pair}.
\end{defi}

\begin{defi}\label{ggdefi}{\rm (\cite{GG})}
\rm
(a) A mereotopological pair $(X,B)$ is called a {\em mereotopological space} if $B$ is a closed base for $X$. We say that $(X,B)$ is a mereotopological $T_0$-space if
$(X,B)$ is a mereotopological space and $X$ is a $T_0$-space.

\smallskip

\noindent(b) A mereotopological space $(X,B)$ is said to be {\em mereocompact} if for every clan $\GA$ of $(B,C_X\cap B^2)$ there exists a point $x$ of $X$ such that $\GA=\s_x^B$ (see (\ref{sg}) for $\sigma_x^B$).
\end{defi}

\begin{rem}\label{remgg}
\rm
(a) Obviously, if $(X,B)$ is a mereotopological space then $X$ is  semiregular.

\smallskip

\noindent(b) The definition of mereocompactness given here is one of the several equivalent expressions of the definition of mereocompactness given in \cite{GG}.  We will use only the definition given here.

\smallskip

\noindent(c) Clearly, a space $X$ is C-semiregular (see Definition \ref{csemi}) iff $(X,RC(X))$ is a mereocompact $T_0$-space. So that, the notion of mereocompactness is an analogue of the notion of C-semiregular space.
\end{rem}

Having in mind the definition of an u-point of a topological space (see Definition \ref{upoint}), we will now introduce the more general  notion of an {\em u-point of a mereotopological pair}.

\begin{defi}\label{mtupoint}
\rm
Let $(X,B)$ be a mereotopological pair and $x\in X$. Then the point $x$ is said to be an {\em u-point of the mereotopological pair $(X,B)$} if, for every $F,G\in B$,
$x\in F\cap G$ implies that $x\in\cl_X(\int_X(F\cap G))$.
\end{defi}

Obviously, a point $x$ of a topological space $X$ is an u-point of $X$ iff it is an u-point of the mereotopological pair $(X,RC(X))$. Also, if $X$ is a topological space, then every point of $X$ is an u-point of the mereotopological pair $(X,CO(X))$.

We will now generalize some assertions from \cite{DV11} concerning u-points.

\begin{pro}\label{remupointme}
 Let $(X,X_0)$ be a topological pair, $(X,B)$ be a mereotopological pair (resp., a mereotopological space) and $B_0=\{F\cap X_0\st F\in B\}$. Then $(X_0,B_0)$ is
 a  mereotopological pair (resp., a mereotopological space); also, for every $x\in X_0$, we have that $x$ is an u-point of $(X_0,B_0)$ iff $x$ is an u-point of $(X,B)$.
\end{pro}

\doc Using Lemma \ref{isombool}, we get that $(X_0,B_0)$ is
 a  mereotopological pair (resp., a mereotopological space). Note that, by Lemma \ref{isombool}, for every  $H\in B$ there exist $H_0\in B_0$ such that $H=\cl_X(H_0)$ and, also, for every  $H_0\in B_0$, $\cl_X(H_0)\in B$. Let $x\in X_0$ be an u-point of $(X_0,B_0)$, $F,G\in B_0$ and $x\in\cl_X(F)\cap\cl_X(G)$. Then $x\in F\cap G$ and hence $x\in\cl_{X_0}(\int_{X_0}(F\cap G))$. Since, by Lemma \ref{isombool}, $\cl_X(cl_{X_0}(\int_{X_0}(F\cap G)))=\cl_X(\int_X(\cl_X(F)\cap\cl_X(G)))$, we get that $x$ is an u-point of $(X,B)$. Conversely, let $x\in X_0$ be an u-point of $(X,B)$. Let $F,G\in B_0$ and $x\in F\cap G$. Then $x\in\cl_X(F)\cap\cl_X(G)$ and thus $x\in\cl_X(\int_X(\cl_X(F)\cap\cl_X(G)))$. Since, by Lemma \ref{isombool}, $X\cap\cl_X(\int_X(\cl_X(F)\cap\cl_X(G)))=\cl_{X_0}(\int_{X_0}(F\cap G))$, we get that $x$ is an u-point of $(X_0,B_0)$.
 \sqs

 \begin{pro}\label{proupointme}
 Let $(X,B)$ be a topological pair and $x\in X$. Then $x$ is an u-point of $(X,B)$ iff $\s_x^B$ is an ultrafilter of the Boolean algebra $B$.
 \end{pro}

 \doc Since, by (\ref{sxb}), $\s_x^B$ is a clan of $(B,C_X\cap B^2)$, we obtain that $\s_x^B$ is a grill of $B$. Hence: ($\s_x^B$ is an ultrafilter of $B$) $\iff$ $[(\fa F,G\in\s_x^B)(F.G\in\s_x^B)]\iff$ ($x$ is an u-point of $(X,B)$). (Note that $F.G=cl_X(\int_X(F\cap G))$.)
 \sqs

 \begin{theorem}\label{csemithme}
 For every  mereocompact $T_0$-space $(X,B)$, the set $$u(X,B)=\{x\in X\st x \mbox{ is an u-point of }(X,B)\}$$ endowed with its subspace topology is a dense zero-dimensional compact Hausdorff subspace of $X$ and is the unique dense zero-dimensional  compact Hausdorff subspace of $X$ such that $RC(X,u(X,B))=B$. Also, the pair $(X,u(X,B))$ is a 2-contact space.
 \end{theorem}

 \doc Set $X_0=u(X,B)$. We have that $\unlb=(B,C_X\cap B^2)$ is a contact algebra. Hence, by Definition \ref{contactspace}(c), $\unlxc(\unlb)$ is a 2-contact space. Let $\unlxc(\unlb)=(X\ap,X_0\ap)$. Then $X\ap=Clans(\unlb)$ and $X_0\ap=Ult(B)$. Since $X$ is a $T_0$-space and $B$ is a closed base for $X$, arguing as in the paragraph after (\ref{txcnc})
 and using (\ref{sxb}),
 we obtain
  that the map
  $$t_X^c:X\lra X\ap,\ \ x\mapsto\s_x^B,$$
is an injection. Also, the fact that $(X,B)$ is mereocompact implies that $t_X^c$ is a surjection.  Now, arguing as in the paragraph immediately after (\ref{txcontc}), we get that
$t_X^c$ is a homeomorphism and
\begin{equation}\label{devsec}
\fa F\in B, \ t_X^c(F)=g_{\unlb}(F).
\end{equation}
 Using Proposition \ref{proupointme}, we obtain that $t_X^c(X_0)=X_0\ap$.
   All this shows that $(X,X_0)$ is a 2-contact space. Hence, $X_0$ is a Stone space and $X_0$ is dense in $X$. Using (\ref{devsec}) and Theorem \ref{mth}(a) (or \cite[Theorem 7.9(a)]{DV11}), we get that $B=RC(X,X_0)=\{\cl_X(P)\st P\in CO(X_0)\}$.

For proving the uniqueness of $X_0$, let $X_0^1$ be a dense Stone subspace of $X$ such that $B=\{\cl_X(P)\st P\in CO(X_0^1)\}$. We will show that $X_0^1=X_0$. Set $B_0=CO(X_0^1)$. Then $B_0=\{F\cap X_0^1\st F\in B\}$ (by Lemma \ref{isombool}) and $(X_0^1,B_0)$ is a mereotopological space. Obviously, every point of $X_0^1$ is an u-point of $(X_0^1,B_0)$. Then Proposition \ref{remupointme} implies that every point of $X_0^1$ is an u-point of $(X,B)$. Thus $X_0^1\sbe X_0$. Obviously, $X_0^1$ is a dense subspace of $X_0$. Since $X_0^1$ is compact and $X_0$ is Hausdorff, we get that $X_0^1=X_0$.
\sqs

Note that Theorem \ref{csemithme} and Proposition \ref{cormthccc} imply the following result from \cite{GG}: if $(X,B)$ is a mereocompact $T_0$-space then $X$ is compact.
(In fact, a stronger result is proved in \cite{GG}: every mereocompact space is compact.)

 \begin{lm}\label{csem2contme}
If $(X,X_0)$ is a 2-contact space,  then $(X,RC(X,X_0))$ is a mereocompact $T_0$-space.
\end{lm}

\doc Clearly, the map $e_{X_0,X}:(CO(X_0),\d_{(X,X_0)})\lra (RC(X,X_0),\ \ F\mapsto\cl_X(F),$ is a CA-isomorphism (use Lemma \ref{isombool}) and for every $x\in X$,
$e_{X_0,X}(\GA_{x,X_0})=\s_x^B$. Then the axioms (CS1), (CS3) and (CS4) from Definition \ref{contactspace} imply that $(X,RC(X,X_0))$ is a mereocompact $T_0$-space.
\sqs

\begin{defi}\label{mcscatdef}
\rm
Let us denote by $\MCS$ the category whose objects are all mereocompact $T_0$-spaces and whose morphisms are all continuous maps between mereocompact $T_0$-spaces which preserve the corresponding u-points (i.e., $$f\in\MCS((X,A),(Y,B))$$ iff $f:X\lra Y$ is a continuous map and, for every u-point $x$ of $(X,A)$, $f(x)$ is an u-point  of $(Y,B)$).
\end{defi}

Obviously, $\MCS$ is indeed a category.

 \begin{theorem}\label{mcsdualthe}{\sc (A Duality Theorem for Contact Algebras)}
 The categories $\CA$ and $\MCS$ are dually equivalent.
 \end{theorem}

 \doc Having in mind Corollary \ref{maindualitytheoremcor2}, it is enough to show that the categories $\CS$ and $\MCS$ are isomorphic. Let $F^g:\CS\lra\MCS$ be defined by $F^g(X,X_0)=(X,RC(X,X_0))$ on the objects of the category $\CS$, and by $F^g(f)=f$ on the morphisms of the category $\CS$. Then Lemma \ref{csem2contme} shows that $F^g$ is well defined on the objects of the category $\CS$. Let $f\in\CS((X,X_0),(Y,Y_0))$. Then $f:X\lra Y$ is a continuous function and $f(X_0)\sbe Y_0$. Since $(X,X_0)$ is a 2-contact space, we get that $X_0$ is a dense Stone subspace of $X$. Since $(X,RC(X,X_0))$ is a mereocompact $T_0$-space and $RC(X,X_0)=\{\cl_X(P)\st P\in CO(X_0)\}$, Theorem \ref{csemithme} implies that $X_0=\{x\in X\st x$  is an u-point of $(X,RC(X,X_0))\}$. Analogously, we get that $Y_0=\{y\in Y\st y$  is an u-point of $(Y,RC(Y,Y_0))\}$. Since $f(X_0)\sbe Y_0$, we get that $f$ preserves the u-points. Hence, $F^g(f)\in\MCS(F^g(X,X_0),F^g(Y,Y_0))$. So, $F^g$ is well defined on the morphisms as well. Obviously, $F^g$ is a functor.

 Let $F^h:\MCS\lra\CS$ be defined by $F^h(X,B)=(X,u(X,B))$ (see Theorem \ref{csemithme} for the notation $u(X,B)$) on the objects of the category $\MCS$, and by $F^h(f)=f$ on the morphisms of the category $\MCS$. Then Theorem \ref{csemithme} implies that $F^h$ is well defined. Obviously, $F^h$ is a functor.

Using once more Theorem \ref{csemithme}, we get that $F^g\circ F^h=Id_{\MCS}$ and $F^h\circ F^g=Id_{\CS}$. Hence, the categories  $\CS$ and $\MCS$ are isomorphic.
\sqs

Finally, we will  show how our results imply the duality for contact algebras described in \cite{GG}.

\begin{defi}\label{gmcscatdef}{\rm (\cite{GG})}
\rm
Let  $\GMCS$ be the category whose objects are all mereocompact $T_0$-spaces and whose morphisms are defined as follows:  $$f\in\GMCS((X,A),(Y,B))$$ iff $f:X\lra Y$ is a  function such that the function $$\psi_f:B\lra A,\ \ F\mapsto f\inv(F),$$ is well defined and is a Boolean homomorphism.
\end{defi}

\begin{cor}\label{gmcsdualthe}{\rm (\cite{GG})}
 The categories $\CA$ and $\GMCS$ are dually equivalent.
 \end{cor}

 \doc We will derive this result from Corollary \ref{maindualitytheoremcor2} showing that the categories $\CS$ and $\GMCS$ are isomorphic. Let $F^i:\CS\lra\GMCS$ be defined by $F^i(X,X_0)=(X,RC(X,X_0))$ on the objects of the category $\CS$, and by $F^i(f)=f$ on the morphisms of the category $\CS$. Then Lemma \ref{csem2contme} shows that $F^i$ is well defined on the objects of the category $\CS$. Let $f\in\CS((X,X_0),(Y,Y_0))$. Then $f:X\lra Y$ is a continuous function and $f(X_0)\sbe Y_0$. Now Corollary \ref{corgtnew1}, (\ref{deftetafc}) and (\ref{deftetafcnov}) show that the function $\psi_f: RC(Y,Y_0)\lra RC(X,X_0),\ H\mapsto f\inv(H),$ is a Boolean homomorphism. Indeed, we have that $H\in RC(Y,Y_0)\iff H=\cl_Y(F)$, where $F\in CO(Y_0)$; thus, for every  $H\in RC(Y,Y_0)$, $H\cap Y_0\in CO(Y_0)$. Therefore, $f\inv(H)=\cl_X(X_0\cap f\inv(H\cap Y_0))$.  Then, using  formula (\ref{deftetafc}) and arguing as in the  paragraph after it, we get that $\psi_f$ is a Boolean homomorphism. Hence $F^i$ is well defined on the morphisms of the category $\CS$ as well. Obviously, $F^i$ is a functor.

 Let $F^j:\GMCS\lra\CS$ be defined by $F^j(X,B)=(X,u(X,B))$ (see Theorem \ref{csemithme} for the notation $u(X,B)$) on the objects of the category $\GMCS$, and by $F^j(f)=f$ on the morphisms of the category $\GMCS$. Then Theorem \ref{csemithme} implies that $F^j$ is well defined on the objects of the category $\GMCS$.
 If $f\in \GMCS((X,A),(Y,B))$ then $f$ is a continuous map because $B$ is a closed base of $Y$ and $f\inv(B)\sbe A$. Let $F^j(X,A)=(X,X_0)$ and $F^j(Y,B)=(Y,Y_0)$. We have to show that $f(X_0)\sbe Y_0$, i.e. that if $x$ is an u-point of $(X,A)$ then $f(x)$ is an u-point of $(Y,B)$. So, let $x$ be an u-point of $(X,A)$, $F,G\in B$ and $f(x)\in F\cap G$. Then $x\in f\inv(F)\cap f\inv(G)$ and $f\inv(F), f\inv(G)\in A$. Hence $x\in f\inv(F). f\inv(G)=\psi_f(F).\psi_f(G)=\psi_f(F. G)=f\inv(F. G)$. Therefore, $f(x)\in F. G$ and thus $f(x)$ is an u-point of $(Y,B)$. So, $F^j$ is well defined on the morphisms of the category $\GMCS$.
 Obviously, $F^j$ is a functor.

 Using once more Theorem \ref{csemithme}, we get that $F^i\circ F^j=Id_{\GMCS}$ and $F^j\circ F^i=Id_{\CS}$. Hence, the categories  $\CS$ and $\GMCS$ are isomorphic.
\sqs

Since, by Theorem \ref{mth}(a), a precontact algebra $\unlb$ is connected iff its canonical 2-precontact space $\unlx(\unlb)$ is connected, all our results from this section which concern dualities have as corollaries duality theorems for the corresponding full subcategories of the corresponding categories of (pre)contact algebras
having as objects all connected (pre)contact algebras. We will now formulate, as an example, the connected variant of Theorem \ref{maindualitytheorem} and will left to the reader the formulation of all other corollaries.

\begin{theorem}\label{mcsdualtheconnnn}
 The full subcategory $\PCAC$ of the category  $\PCA$ whose objects are all connected precontact algebras  is dually equivalent to the full subcategory $\PCSC$ of the category $\PCS$ having as objects all connected 2-precontact spaces.
 \end{theorem}


\begin{thebibliography}{99}

\parskip0pt

\bibitem{AHS}
 {\sc   Ad\'amek, J., Herrlich, H. and  Strecker, G. E.}
  \newblock  {\em Abstract and Concrete Categories},

  \bibitem{AP}
{\sc    Alexandroff, P. S. and Ponomarev, V. I.}
\newblock  On bicompact extensions of topological spaces,
\newblock  Vestn. Mosk. Univ. Ser. Mat. (1959), 93-108. (In Russian).



\bibitem{BBSV}
{\sc Bezhanishvili, G., Bezhanishvili, N., Sourabh, S. and Venema, Y.}
\newblock Subordinations, closed relations, and compact Hausdorff spaces.
\newblock (December 2014) pp. 1-21, available  at

http://www.phil.uu.nl/\~{}bezhanishvili/Papers/stone\_closed\_relation\_9.pdf

\bibitem{CE}
{\sc  \v Cech, E.}
\newblock  {\em Topological Spaces.}
\newblock  Interscience, London, 1966.


\bibitem{CNG}
{\sc Comfort, W.  and   Negrepontis, S.}
\newblock   {\em Chain Conditions in Topology},
\newblock Cambridge Univ. Press, Cambridge, 1982.

\bibitem{deV}
 {\sc de Vries, H.}
\newblock  {\em Compact Spaces and
Compactifications, an Algebraic Approach},
\newblock Van Gorcum, The Netherlands, 1962.


\bibitem{D-APCS09}
{\sc Dimov, G.}
\newblock  	A generalization of De Vries' Duality Theorem,
\newblock  Applied Categorical Structures 17 (2009), 501-516.


\bibitem{D2009}
{\sc Dimov, G.}
\newblock Some generalizations of the Fedorchuk duality theorem - {I},
\newblock Topology Appl., 156 (2009), 728-746.

\bibitem{D-AMH1-10}
{\sc Dimov, G.}
\newblock  A de {V}ries-type duality theorem for the category of locally compact
  spaces and continuous maps - {I},
  \newblock  Acta Math. Hungarica, 129 (4) (2010), 314-349.

  \bibitem{D-AMH2-11}
{\sc Dimov, G.}
\newblock  A de {V}ries-type duality theorem for the category of locally compact
  spaces and continuous maps - {II},
  \newblock  Acta Math. Hungarica, 130 (1) (2011), 50-77.


\bibitem{D2012}
{\sc Dimov, G.}
\newblock Some Generalizations of the Stone Duality Theorem,
\newblock Publicationes Mathematicae Debrecen 80(3-4) (2012), 255-293.


\bibitem{DI2016}
{\sc Dimov, G. and E. Ivanova}
\newblock Yet another duality theorem for locally compact spaces,
\newblock Houston Journal of Mathematics (2016) (to appear).

\bibitem{DV1}
{\sc Dimov, G. and Vakarelov, D.}
\newblock Contact Algebras and Region-based Theory of Space: A
Proximity Approach - I. Fundamenta Informaticae 74(2-3) (2006), 209-249.

\bibitem{DV2}
{\sc Dimov, G. and Vakarelov, D.}
\newblock Contact Algebras and Region-based Theory of Space: A
Proximity Approach - II. Fundamenta Informaticae  74 (2-3) (2006), 251-282.

\bibitem{DV3}
{\sc Dimov, G. and Vakarelov, D.}
\newblock Topological Representation of Precontact Algebras.
\newblock In: {\em Relation Methods in Computer Science}, W. MacCaull, M. Winter, I. Duentsch (Eds.),
Lecture Notes in Computer Science, 3929 (2006), 1-16, Springer-Verlag Berlin Heidelberg.

\bibitem{DV11}
{\sc Dimov, G. and Vakarelov, D.}
\newblock Topological Representation of Precontact Algebras and a Connected Version of the Stone Duality Theorem -- I.
\newblock arXiv:1508.02220v3, 1-44 (a revised version of this paper will appear in Topology Appl.).


\bibitem{DUV}
{\sc D{\"u}ntsch, I. and Vakarelov, D.}
\newblock Region-based theory of discrete spaces: A proximity
approach.
\newblock In: Nadif, M., Napoli, A., SanJuan, E., and Sigayret, A.
EDS,
\newblock {\em Proceedings of Fourth International Conference
Journ{\'e}es de l'informatique Messine}, 123-129,
\newblock Metz, France, 2003.
\newblock Journal version in: {\em Annals of Mathematics and
Artificial Intelligence},  49(1-4) (2007), 5-14.


\bibitem{DW}
{\sc D{\"u}ntsch, I. and Winter, M.}
\newblock A Representation theorem for Boolean Contact Algebras.
\newblock Theoretical Computer Science (B), 347 (2005), 498-512.

\bibitem{EF}
{\sc Efremovi\v{c},  V. A.}
\newblock Infinitesimal spaces.
\newblock  DAN SSSR, 76 (1951), 341--343.


\bibitem{E}
{\sc Engelking, R.}
\newblock {\em  General Topology},
\newblock PWN, Warszawa, 1977.

\bibitem{F}
 {\sc Fedorchuk,  V. V.}
\newblock   Boolean $\d$-algebras and quasi-open mappings.
\newblock  Sibirsk. Mat. \v{Z}. 14 (5) (1973), 1088--1099; English translation: Siberian Math. J.
14 (1973), 759-767 (1974).



\bibitem{Galton}
{\sc Galton, A.}
\newblock The mereotopology of discrete spaces.
\newblock In: Freksa, C. and Mark, D.M. eds,
\newblock {\em Spatial Information Theory, Proceedings of the
International Conference COSIT'99,}
\newblock Lecture Notes in Computer Science, 251--266,
\newblock Springer-Verlag, 1999.

\bibitem{GG}
{\sc  Goldblatt, R.  and  Grice, M.}
\newblock Mereocompactness and duality for
mereotopological spaces.
\newblock In:  Katalin Bimbo (Ed.), J. Michael Dunn on Information Based Logics,  Springer, 2016 (to appear).

Available at
http://homepages.ecs.vuw.ac.nz/~rob/papers/mereo.pdf


\bibitem{kop89}
{\sc  Koppelberg, S.}
\newblock  {\em Handbook
  on Boolean Algebras, vol. 1: General Theory of Boolean Algebras},
\newblock  North Holland,
1989.

\bibitem{NW}
Naimpally, S.  and Warrack, B.
{\it Proximity Spaces.}
Cambridge,
London, 1970.





\bibitem{Sm}
{\sc Smirnov, J. M. }
\newblock On proximity spaces.
\newblock Mat. Sb. 31 (1952), 543--574.

\bibitem{ST}
{\sc Stone, M. H.}
\newblock The theory of representations for Boolean algebras.
\newblock    Trans. Amer. Math. Soc.,  40  (1936), 37--111.



\bibitem{Thron}
{\sc Thron, W.}
\newblock Proximity structures and grills.
\newblock {\em Math. Ann.}, 206 (1973), 35-62.

\bibitem{VDDB}
{\sc Vakarelov, D., Dimov, G., D{\"u}ntsch, I., Bennett, B.}
\newblock A proximity approach to some region-based theory of
space.
\newblock {\em Journal of applied non-classical logics},  12
 (3-4) (2002), 527-559.





\end{thebibliography}
\end{document}